\documentclass[12pt]{amsart}

\usepackage{amssymb,amsmath,amsthm}  
\usepackage{graphicx}
\usepackage{bm}
\usepackage{enumerate}
\usepackage{enumitem}
\usepackage{braket}
\usepackage{amsaddr}
\usepackage{amscd}
\usepackage[dvipsnames]{xcolor}
\usepackage[mathscr]{eucal}
\usepackage{epstopdf}
\usepackage[hidelinks]{hyperref}
\usepackage{subcaption}
\usepackage{cleveref}
\usepackage{relsize}
\hypersetup{
  colorlinks   = true, 
  urlcolor     = PineGreen, 
  linkcolor    = PineGreen, 
  citecolor   = PineGreen 
}

 \allowdisplaybreaks
   \DeclareMathOperator*{\argmin}{argmin}
\numberwithin{equation}{section}
\newtheorem{theorem}{Theorem}[section]

\newtheorem{assumption}[theorem]{Assumption}

\newcommand{\tu}{\textup}

\newcommand{\dd}{\displaystyle}
\newcommand{\bfa}[1]{\boldsymbol{#1}} 			%

\usepackage[numbers]{natbib}
\usepackage[latin9]{inputenc}
\usepackage[letterpaper]{geometry}
\usepackage{pifont}
\usepackage{float}
\usepackage{color}
\usepackage{adjustbox}
\usepackage{diagbox}
\usepackage{multirow}
\usepackage[makeroom]{cancel}
\usepackage{array}
\usepackage{mathrsfs}
\usepackage{tikz}
\usepackage{comment}
\usepackage{csquotes}

\graphicspath{{./Figures/}}

\usepackage{amsfonts}

\usepackage{color}
\usepackage{adjustbox}
\usepackage{diagbox}
\definecolor{black}{rgb}{0,0,0}

\definecolor{red}{rgb}{1,0,0}

\definecolor{blue}{rgb}{0,0,1}

\newcommand{\todo}[1]{{\color{red}{#1}}}

\numberwithin{equation}{section}

\renewcommand{\div}{\mathop{\rm div}\nolimits}

\newcommand{\dx}{ \mathrm{d}x}

\newcommand{\beq}{\begin{equation}}
\newcommand{\eeq}{\end{equation}}
\newcommand{\beqq}{\begin{equation*}}
\newcommand{\eeqq}{\end{equation*}}
\newcommand{\beqas}{\begin{eqnarray*}}
\newcommand{\eeqas}{\end{eqnarray*}}
\newcommand{\bsp}{\begin{split}}
\newcommand{\eesp}{\end{split}}

\makeatother

\usepackage[english]{babel}

\date{\today}

\title[P\MakeLowercase{rediction of} \MakeLowercase{numerical homogenization} \MakeLowercase{using} \MakeLowercase{deep} \MakeLowercase{learning} \MakeLowercase{for the} R\MakeLowercase{ichards equation}]
{\textsf{\LARGE P\MakeLowercase{rediction of} \MakeLowercase{numerical homogenization} \MakeLowercase{using} \MakeLowercase{deep} \MakeLowercase{learning} \MakeLowercase{for the} R\MakeLowercase{ichards equation}}}

\author[S\MakeLowercase{ergei} S\MakeLowercase{tepanov}, D\MakeLowercase{enis} S\MakeLowercase{piridonov}, T\MakeLowercase{ina} M\MakeLowercase{ai}]{S\MakeLowercase{ergei} S\MakeLowercase{tepanov$^a$}, D\MakeLowercase{enis} S\MakeLowercase{piridonov$^a$}, T\MakeLowercase{ina} M\MakeLowercase{ai$^{b,c,*}$}
}

\begin{document}


\begin{abstract}
For the nonlinear Richards equation as an unsaturated flow through heterogeneous media, we build a new coarse-scale approximation algorithm utilizing numerical homogenization. 
This approach follows deep neural networks (DNNs) to quickly and frequently calculate macroscopic parameters. 
More specifically, we train neural networks with a training set consisting of stochastic permeability realizations and corresponding computed macroscopic targets (effective permeability tensor, homogenized stiffness matrix, and right-hand side vector). Our proposed deep learning scheme develops nonlinear maps between such permeability fields and macroscopic characteristics, and the treatment for Richards equation's nonlinearity is included in the predicted coarse-scale homogenized stiffness matrix, which is a novelty.  This strategy's good performance is demonstrated by several
numerical tests in two-dimensional model problems, for predictions of the macroscopic properties and consequently solutions.
\end{abstract}

\maketitle


\vskip10pt

%



\vfill

\noindent \textit{Sergei Stepanov}; $^a$Laboratory of Computational Technologies for Modeling Multiphysical and Multiscale Permafrost Processes, North-Eastern
Federal University, 677000 Yakutsk, Republic of Sakha (Yakutia), Russia;  \texttt{cepe2a@inbox.ru}\\

\noindent \textit{Denis Spiridonov}; $^a$Laboratory of Computational Technologies for Modeling Multiphysical and Multiscale Permafrost Processes, North-Eastern
Federal University, 677000 Yakutsk, Republic of Sakha (Yakutia), Russia; \texttt{d.stalnov@mail.ru}\\



\noindent $^*$Corresponding author: \textit{Tina Mai}; $^b$Institute of Research and Development, Duy Tan University, Da Nang, 550000, Vietnam; $^c$Faculty of Natural Sciences, Duy Tan University, Da Nang, 550000, Vietnam; \texttt{maitina@duytan.edu.vn}


\newpage













\section{Introduction}\label{intro}

Unsaturated flow happens in soils when soil water travels downward, the topsoil dries out due to evaporation, or water permeates the soil surface across many soils' pores \cite{unflow, soilmoisture}.
This generally leads to the conception of soil moisture \cite{soilmoisture}, described by water content (the amount of water in unsaturated soils) and water potential (the water's energy state).
Despite being a minor component of the water cycle, soil moisture is necessary for a number of hydrological, biological, and biogeochemical processes.  For instance, soil moisture is a crucial factor in agronomy, environmental management, sanitary engineering, groundwater storage, weather prediction, balances of energy, earth's dynamic systems, etc.  The related procedures are mathematically modeled by the Richards equation \cite{r0,modifiedPicard,richarde2,haverkamp,richardsreview},
depicting how water moves through unsaturated soil surface that is filled with both water and air \cite{ry1, soilmoisture}.

In this paper, the Richards equation is considered in non-periodic and heterogeneous media.  Obstacles further develop when there are material property uncertainties in some local areas, which are typical for subsurface activities, oil reservoirs, infiltration, aquifers, and most real-world simulation models. The variability of such properties' uncertainty is usually spatial. 

Traditional direct methods can solve those problems in a few steps using fine-grid simulation.  First, a local fine grid is built.
The equation is then discretized on that fine grid, and the set of local solutions is utilized to form a global solution.  This approach can be operated with well-known frameworks like the Finite Element Method (FEM) \cite{fem} and the Finite Volume Method (FVM).  The grid size needs to be as small as possible to explicitly resolve all grid-level heterogeneity in order for it to be applicable and to achieve convergence.  However, the smaller the grid size is, the larger the computing resources are demanded.  Hence, some sort of model reduction is required. Familiar techniques involve dividing the focused domain into coarse-scale grid blocks and then calculating the effective coefficients in each coarse block \cite{homod}. This computation needs the local problems' fine-scale solutions in each coarse-grid cell or representative volume (in conventional upscaling approaches based on homogenization, whose expansion can be used to derive multicontinuum methods in a new article \cite{yewl22homo}).  For our case of non-periodic media, numerical homogenization \cite{homo1} is applied to single-continuum Richards equation.
Concerning approximation on a macroscale computational grid, homogenization assumes that local problems' solutions carry microscale information. 

Given largely different permeability fields' realizations, the computing effort can become enormous.  Therefore, establishing a functional link between the permeability fields and the local macroscopic coefficients can prevent the need for repeated, costly calculations and so significantly alleviate the computational complexity.  Because of various media properties, such a functional relationship is nonlinear and thus usually needs high-order approximations when modelling.  Consequently, utilizing machine learning techniques to create these complicated models is natural \cite{dnngms19,ex-imML}.

Deep neural networks (DNNs) fall under the machine learning framework \cite{dnneg14}, which constitutes a type of artificial intelligence. In this research, we use DNNs to investigate numerical homogenization.  A deep neural network's structure is often made up of numerous layers, each of which contains several neurons (as nonlinear processing units for extracting features) \cite{dnnnature}. Right after the input layer, the first hidden layer of a deep neural network (DNN) learns basic characteristics, which are then supplied to the second hidden layer, for training itself to identify increasingly complex and abstract properties. This architecture continues on succeeding levels as the layers build up until the output is accurate enough to be considered acceptable.
%
To enhance the expressive capability of neural networks, some nonlinear activation functions (such as ReLU, sigmoid, and tanh) are required between layers, from an input neuron to an output neuron. The output is then sent into the network's next layer as an input. Recently, deep neural networks have been effectively applied to problems needing pattern recognition, including picture identification, voice recognition as well as processing of natural language, and have been used to comprehend complex data sets \cite{dnneg8,dnneg9,dnneg10}.  Additionally, extensive study has been carried out to examine the descriptive potential of deep neural networks \cite{dnneg11,dnneg12,dnneg13,dnneg14,dnneg15}.  There is a variety of other applications employed DNN \cite{nlnlmc31,wang2020deep,wang2020reduced, dnneg20} and convolutional neural network (CNN) \cite{cnn-homo-elasticity19, cnn-homo-poro20}, which also motivate this work. 

In our paper, neural networks are utilized to predict macroscopic parameters (effective permeability, coarse-grid homogenized stiffness matrix, and right-hand side vector). Here, a machine learning method is developed using deep neural network and graphics processing unit (GPU) during the training procedure via two experiments. 
In particular, we predict an effective permeability in the first experiment and predict a coarse-grid homogenized matrix together with a right-hand side vector in the second experiment,  where the predicted homogenized stiffness matrix contains information on how the Richards equation's nonlinearity is tackled, which is a new contribution. For various input random permeability fields, a corresponding set of targets (desired output) is employed to train the neural networks. Fast and accurate computations can be conducted once neural networks have been trained on the training data. Our numerical results show that the deep networks produce favorable results and that they can be applied well to the testing data. 

The paper is organized as follows.   In Section \ref{fus}, we begin with Subsection \ref{sec:model} about model problem for the underlying Richards equation that characterizes the unsaturated flow within non-periodic heterogeneous media, then Subsection \ref{pre} regards coarse-scale discretization as well as Picard iteration for linearization, and the numerical homogenization for nonlinear Richards equation is introduced in Subsection \ref{sec:method}.  Next, the concept of employing deep learning as a stand-in for numerical homogenization prediction is described in Section \ref{dlhomo}, where the sampling is thoroughly discussed, and the networks are carefully defined.  In Section \ref{numer}, we give numerical examples (especially for solutions) to demonstrate the performance of our proposed networks, using predictions of effective permeability and coarse-grid homogenized stiffness matrix along with right-hand side vector, by several cases with various configurations.  Section \ref{sec:conclude} provides a summary and conclusion.  Global convergence of Picard linearization process is shown in Appendix \ref{cp}.

\section{Preliminaries}\label{fus}


\subsection{Model problem}\label{sec:model}

In $\mathbb{R}^d\,,$ let $\Omega$ be a computational domain that is open, bounded, convex, Lipschitz, and simply connected.
Throughout this work, the case $d = 2$ is investigated to make our discussion clearer, but the concept can be simply extended to $d = 3$ as well.
Latin indices (as $i,j$) are in $\{1,2\}\,.$  The symbols $\nabla$ and $\dfrac{\partial}{\partial t}$ respectively indicate the spatial gradient and temporal derivative, across this paper.  Further notation follows \cite{mcl,rtt21}.  
Bold letters (e.g., $\bfa{v}$ and $\bfa{T}$) stand for vector fields and matrix fields over $\Omega\,,$ while italic capitals (e.g., $f$) are used to represent functions.  The spaces of functions, vector fields, and matrix fields over $\Omega$ are described by italic capitals (e.g., $L^2(\Omega)$), 
boldface Roman capitals (e.g., $\bfa{V}$), 
and special Roman capitals (e.g., $\mathbb{S}$), respectively. 

The $L^2$ inner product is symbolized by $(\cdot,\cdot)\,$, and the Sobolev space 
$V: = H_0^1(\Omega) = W_0^{1,2}(\Omega)\,$ has the norm $\| \cdot \|_{V}$: 
\[\|v\|_{V} = \left(\|v\|^2_{L^2(\Omega)} + 
\|\nabla v \|^2_{\bfa{L}^2(\Omega)}\right)^{\frac{1}{2}}\,.\]
Herein, $\| \nabla v \|_{\bfa{L}^2(\Omega)}:= \| | \nabla v | \|_{L^2(\Omega)}\,,$ where 
$| \nabla v|$ specifies the Euclidean norm of the $2$-component vector-valued function $ \nabla v\,.$  We also let $\bfa{V} = V^2 = [H_0^1(\Omega)]^2\,.$  For any $\bfa{v} = (v_1,v_2) \in \bfa{V}\,,$ it holds that $\| \nabla \bfa{v}\|_{\mathbb{L}^2(\Omega)}:= \| | \nabla \bfa{v}| \|_{L^2(\Omega)}\,,$ having 
$| \nabla \bfa{v}|$ as the Frobenius norm of the $2 \times 2$ matrix $\nabla \bfa{v}\,.$ 

\bigskip

In the literature, the form of time-dependent Richards equation \cite{r0} is usually 
\beq
\label{eq:originaltheta}
\frac{\partial \Theta(p(t,\bfa{x}))}{\partial t} - \div (\varkappa(\bfa{x},p(t,\bfa{x}))\nabla p(t,\bfa{x}))
	= f(t,\bfa{x}) \  \textrm{in} \ (0,T] \times \Omega\,,
	\eeq
where $T>0$ denotes the terminal time.  This equation has the initial condition $p(0,\bfa{x})= p_{0}$ in $\Omega$ and the Dirichlet boundary condition $p(t,\bfa{x})=0$ on $(0,T] \times \partial \Omega\,.$  The fundamental notation can be found in \cite{rtt21}.
Here, volumetric soil water content is denoted by $\Theta(p(t,\bfa{x}))\,,$ 
the pressure head is described by $p:=p(t,\bfa{x})\,,$ the unsaturated hydraulic conductivity is represented by $\varkappa(p):=\varkappa(\bfa{x},p)\,,$ and the source or sink function is abbreviated by $f\,.$

As in \cite{rtt21}, we remark that the volumetric water content $\Theta(p)$ in \eqref{eq:original0} is often a nonlinear function of the pressure head $p\,,$ and it has the form as follows \cite{modifiedPicard}:
\[\frac{\partial \Theta(p)}{\partial t} = C(p) \frac{\partial p}{\partial t}\,.\]
The nonlinear hydraulic conductivity $\varkappa(\bfa{x}, p)$ is important for homogenization in this paper.  On the other hand, $C(p)$ part is not necessary and can be eliminated, resulting in the identity function $\Theta(p) = p$ in our considering Richards equation:
\beq
\label{eq:original0}
\frac{\partial p(t,\bfa{x})}{\partial t} - \div (\varkappa(\bfa{x},p(t,\bfa{x}))\nabla p(t,\bfa{x}))
= f(t,\bfa{x}) \  \textrm{in} \ (0,T] \times \Omega\,.
\eeq
Later, one can see that Eq.\ \eqref{eq:original0} has a multiscale high-contrast coefficient $\kappa(\bfa{x}) = \kappa(\bfa{x},\omega)\,,$ which is generally stochastic \cite{cnn-mcmc20}.  But for simplicity, the stochasticity notation $\omega$ is omitted.

We assume that the hydraulic conductivity and its spatial gradient are uniformly bounded, particularly, there exist positive constants $\underline{\varkappa}$ and $\overline{\varkappa}$ satisfying the following inequalities: 
\begin{align}
	\label{Coercivity}
	\begin{split}
		\underline{\varkappa} \leq \varkappa(\bfa{x}, p), \ | \nabla \varkappa(\bfa{x},p)| \leq \overline{\varkappa}\,.
	\end{split}
\end{align}
Also, it is assumed that the initial condition is 
\begin{equation}\label{ini}
	p_0 = p(0,\bfa{x}) \in V\,,
\end{equation}
without losing generality.

The following bilinear form is defined with a given $u \in V$: 
for all $p\,, v \in V\,,$
\begin{align}\label{ai}
	a(p,v;u)=\int_{\Omega} \varkappa(u)\nabla p \cdot \nabla v \, \dx\,.
\end{align}
Then, \eqref{eq:original0} posesses the variational form as follows:  
find $p \in V$ such that 
\begin{align}\label{r1e}
	\left(\frac{\partial{p}}{\partial t} , v \right) + a(p,v;p) = (f,v)\,,
\end{align}
with all $v \in V\,,$ for a.e.\ $t \in (0,T]\,,$ and $f(t,\cdot) \in L^2(\Omega)\,.$  In \eqref{ini}, the initial condition was specified.

\subsection{Coarse-scale discretization and Picard iteration for linearization}
\label{pre}

We utilize an efficient Picard iterative technique as described in \cite{rpicardc, Spiridonov2019, cemnlporo,tfcmm,ttr22} to deal with the nonlinearity of our problem. For the time-dependent Richards equation, such an iterative algorithm is provided in this section.



To achieve the first purpose of \eqref{r1e}'s time discretization (see, for example, \cite{rpicardc, Spiridonov2019}), we shall employ the classical backward Euler finite-difference method: find $p \in V$ such that 
for all $v \in V\,,$
\begin{align}\label{r1ed}
	\left(\frac{p_{s+1} - p_{s}}{\tau} , v \right) + a(p_{s+1},v;p_{s+1})= (f_{s+1},v)\,, 
\end{align}
where we separate the time domain $[0,T]$ uniformly into $S$ intervals, the size of temporal step is $\tau = T/S > 0$, 
and a function's evaluation at the time $t_s= s\tau$ is designated by the subscript $s$ (for $s=0,1,\cdots,S$).

The nonlinearity in space will then be tackled by Picard linearization iteration (see, for example, \cite{rpicardc, Spiridonov2019,cemnlporo}) as follows.  Provided $p_s\,,$ at the time step $(s+1)$th, $p^0_{s+1} \in V$ is guessed.  
With $n=0,1,2, \cdots\,,$ find $p^{n+1}_{s+1} \in V$ such that for all $v \in V\,,$ 
\begin{align}\label{r1el}
	\left(\frac{p^{n+1}_{s+1} - p_{s}}{\tau} , v \right) + a(p^{n+1}_{s+1},v;p^n_{s+1}) 
	= (f_{s+1},v)\,. 
\end{align}
Equivalently,
\begin{equation}\label{r1elr}
\left(\frac{p^{n+1}_{s+1}}{\tau}\,, v \right)  + a(p^{n+1}_{s+1},v;p^{n}_{s+1}) = \left(\frac{p_{s}}{\tau}\,, v \right) + (f_{s+1},v) \,.
\end{equation}

In \cite{rh2}, one can find the proof for existence and uniqueness of solution $p_{s+1}^{n+1}$ to the linearized equation \eqref{r1el}. 
%
%
As $n$ comes to $\infty\,,$ the Picard iterative process converges to a limit, which is theoretically proved in Appendix \ref{cp}. 
%
%
Numerically, this process is finished when it approaches some $\alpha$th iteration at a specified stopping requirement.  To go to the next temporal step in \eqref{r1ed}, the prior time data is set as
\begin{equation}\label{pdata}
	p_{s+1} = p_{s+1}^{\alpha}
\end{equation}
All over this paper, we propose a halting indicator based on the relative successive difference, which states that provided a user-defined tolerance $\delta_0 > 0\,,$ if 
\begin{equation}\label{pt}
\dfrac{\|p_{s+1}^{n+1} - p_{s+1}^{n} \|_{L^2(\Omega)}}{\| p_{s+1}^{n} \|_{L^2(\Omega)}} \leq \delta_0\,,
\end{equation}
then the iteration procedure is ended.  In Section 
\ref{numer}, $\delta_0 = 10^{-6}$ is chosen. 

The fine-grid notation is now discussed.  In order to solve the local problem \eqref{epk0}--\eqref{epk}, we let $\mathcal{T}_h$ (\textit{fine grid}) be a conforming triangular partition for the computational domain $\Omega$, having local grid sizes $h_P: = \tu{diam}(P)\; \forall P \in \mathcal{T}_h$, and $h:= \displaystyle \max_{P \in \mathcal{T}_h} h_P$.
The size $h$ is assumed to be small enough so that the fine-grid solution is reasonably near the correct solution.  In numerical simulation, the fine grid possesses equal squares, where we employ a fine-scale triangulation whose two triangle components are built inside each square cell, as shown in Fig. \ref{fig:fine1block} below.  It is worth noting that the fine grid $\mathcal{T}_h$ is only used to numerically handle local problems.
\begin{figure}[h!]
	\begin{center}
	\begin{minipage}[h]{0.27\linewidth}
			\center{\includegraphics[width=\linewidth]{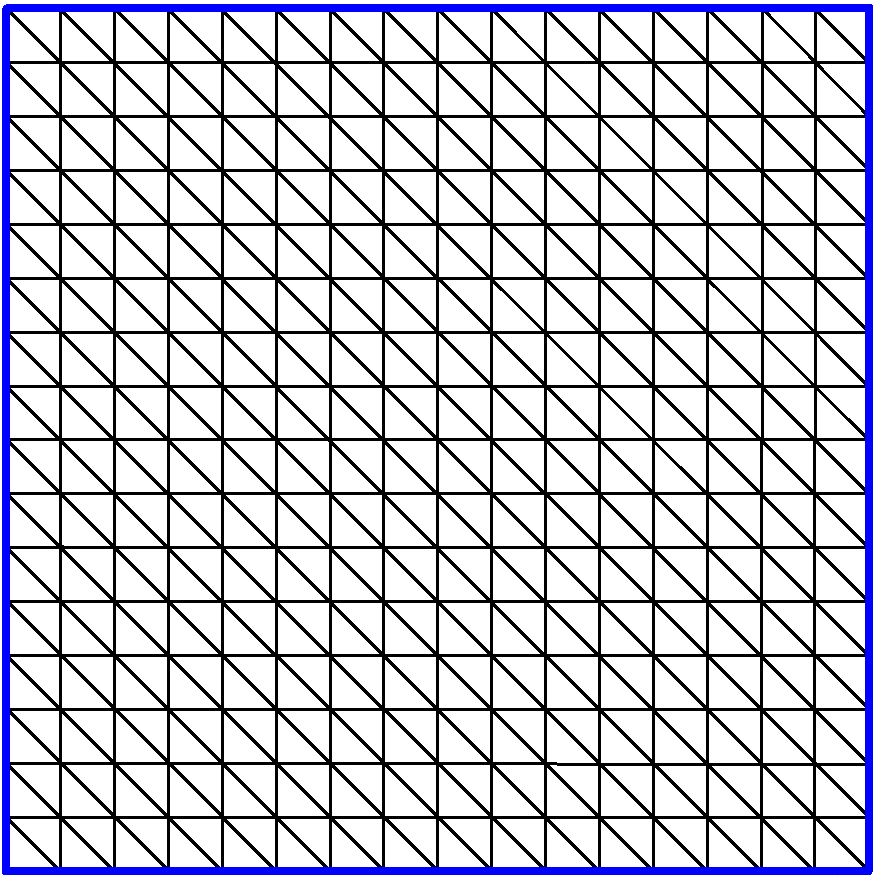}}
		\end{minipage}
	\end{center}
	\caption{Illustration of a triangular fine grid in one coarse block.}
	\label{fig:fine1block}
\end{figure}

\begin{figure}[ht!]
	\centering
	\includegraphics[width=0.30\linewidth]{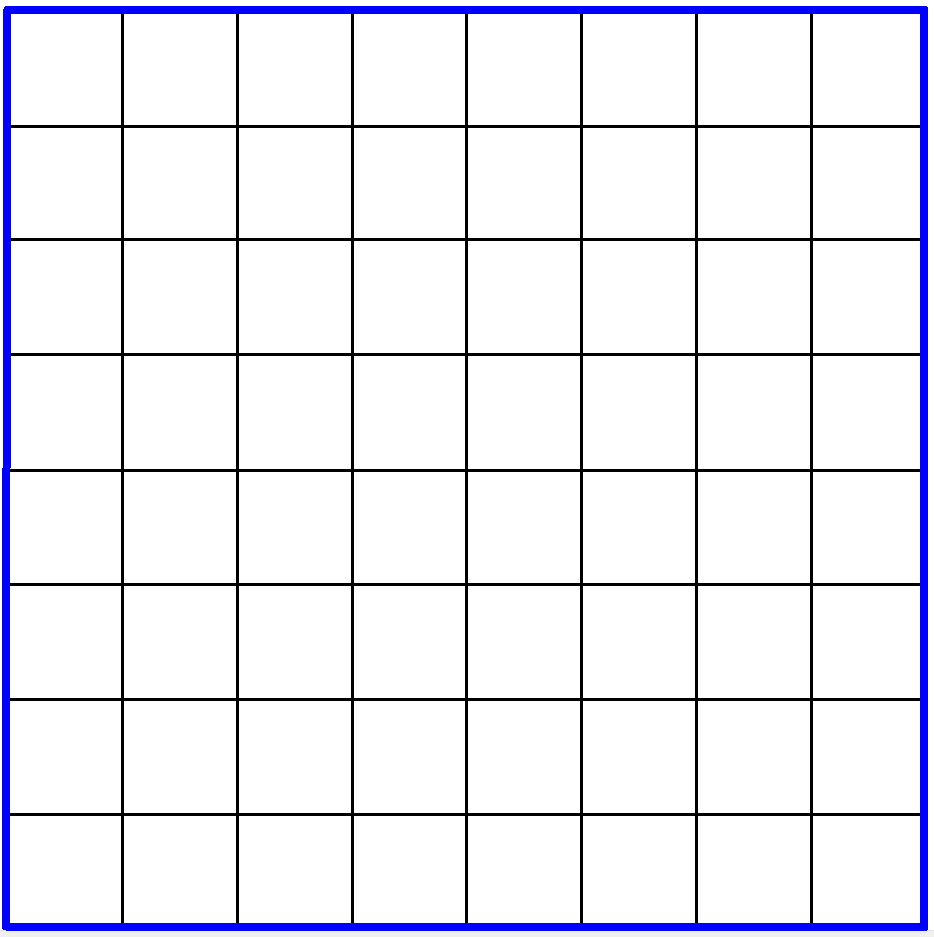}
	\caption{Illustration of a coarse grid, where each smallest square is a coarse block.}
	\label{fig:mesh}
\end{figure}

Next, we review the coarse grid conception, which will be used to discretize the variational problem \eqref{r1e}.  The first step is dividing $\Omega$ into finite elements. 
This partition is known as coarse grid $\mathcal{T}^H\,,$ which has $\mathcal{T}_h$ as a refinement.  An element $K$ in $\mathcal{T}^H$ is referred to as a coarse-grid block (which is also called coarse-grid cell, coarse element, or coarse patch).  Furthermore, the coarse-grid size $H$ is defined as the side length of each square coarse cell,
with $H \gg h\,.$
The number of coarse blocks is $N_b\,,$ while the number of coarse-grid nodes is $N_v\,.$ Then, $\{\bfa{x}_k\}_{k=1}^{N_v}$ exhibits the set of all coarse nodes (vertices).  A coarse grid is shown in Figure~\ref{fig:mesh}.  

\bigskip

With regard to the coarse grid $\mathcal{T}^H\,,$ let $V^H$ be the first-order Galerkin (standard) finite element basis space:
%
\[V^H:= \{ v \in V: v|_P \text{ is a polynomial of degree } \leq 1 \; \forall P \in \mathcal{T}^H\}\,.\]
In this $V^H\,,$ provided an initial $p_{H,0} = P_H p_0\,,$ where $p_0$ is from \eqref{ini} and $P_H$ is the $L^2(\Omega)$ projection operator onto $V^H\,.$  

The fully Picard discrete scheme in $\mathcal{T}^H$ is as follows: We start with an initial $p_{H,0} \in V^H\,,$ and guess $p^0_{H,s+1}\in V^H$ at the $(s+1)$th time step, then do iteration in $V^H$ for $p^{n+1}_{H,s+1}\in V^H\,:$
\begin{align}\label{r1elH}
	\left(\frac{p^{n+1}_{H,s+1} - p_{H,s}}{\tau} , v \right) + a(p^{n+1}_{H,s+1},v;p^n_{s+1})= (f_{H,s+1},v)\,,
\end{align}
with all $v \in V^H$, and $n=0,1,2, \cdots\,,$ until running into \eqref{pt} at some Picard step $\alpha$th.  
To advance to the next time step in \eqref{r1ed}, we employ \eqref{pdata} to establish the prior time data $p_{H,s+1} = p_{H,s+1}^{\alpha}\,.$


\subsection{Numerical homogenization}\label{sec:method}


The purpose of this section is to obtain numerical homogenization of the nonlinear Richards equation \eqref{eq:original0}
as coarse-grid approximation.  We follow the numerical homogenization scheme in \cite{homo1,15mualemv} (nonlinear) and thanks to \cite{homo-poro19} (linear).
More specifically, the effective coefficient is first numerically calculated.  The homogenized Richards equation is then solved on coarse grid as
in Subsection \ref{pre}, where the nonlinearity may be treated as a constant at each Picard iteration step (after temporal discretization).  We denote the volume average over any $K \in \mathcal{T}^H$ as follows:
\[\langle \cdot \rangle_K = \frac{1}{|K|} \int_K (\cdot) \, \dx\,.\]

\subsubsection{Effective coefficient}\label{ep}
The effective coefficient $\tilde{\varkappa}$ (also known as homogenized coefficient or effective parameter \cite{homo1}) 
is computed numerically in this work.  We consider the Richards equation \eqref{eq:original0} in non-periodic heterogeneous media.  The piecewise constant nonlinear function $\tilde{\varkappa}(\bfa{x},s)$ in variable $\bfa{x}$ over $\mathcal{T}^H$ is defined for any $s \in \mathbb{R}$ as follows. 

First, let $\psi_j(\bfa{x})$ ($j=1,\cdots,d$) be the fine-scale solution of the following problem for any $K \in \mathcal{T}^H$ and $s \in \mathbb{R}$:
\begin{align}
-\div (\varkappa(\bfa{x},s) \, \nabla \psi_{j}) &= 0 \ \ \ \text{in} \ K, \label{epk0}\\
\psi_j &= x_j  \ \ \ \text{on} \ \partial K\,, \label{epk}
\end{align}
having $\varkappa(\bfa{x},s)$ as the random heterogeneous coefficient in $\bfa{x}$ restricted to the local region $K \in \mathcal{T}^H\,,$ with $\bfa{x} = (x_1,\cdots, x_d)\,,$ for $d=2\,.$  It is worth noting that this local problem is solved on the fine grid $\mathcal{T}_h\,.$

Then, with $j = 1,\cdots,d\,,$ the constant tensor $\tilde{\varkappa}(s)$ is computed on the current $K$ by 
\begin{align}\label{k8}
\tilde{\varkappa}(s) \langle \nabla \psi_j \rangle_K = \langle \varkappa(\bfa{x},s) \nabla \psi_j \rangle_K \,. 
\end{align}
As a result, for $j,l = 1,\cdots,d\,,$ each entry in the constant matrix $\tilde{\varkappa}(s)$ on $K$ has the form \cite{homo-poro19}
\begin{align}\label{k8e}
	\tilde{\varkappa}_{jl}(s) = \frac{1}{|K|} \int_K \varkappa(\bfa{x},s) \frac{\partial \psi_j}{\partial x_l} \dx\,.
\end{align}

Using the Divergence Theorem
and \eqref{epk}, it holds that \eqref{k8} is well-defined:
\begin{equation}\label{green}
\langle \nabla \psi_j \rangle_K = \frac{1}{|K|} \int_{K} \nabla \psi_j \, \dx = \frac{1}{|K|} \int_{\partial K} \psi_j \, \bfa{n} \,  \tu{d}\sigma =  \bfa{e}_j + \frac{1}{|K|} \int_{\partial K} (\psi_j - x_j) \, \bfa{n} \,  \tu{d}\sigma
=  \bfa{e}_j\,,  
\end{equation}	
where $\bfa{n}$ represents the piecewise smooth field 
(on $\partial K$) of outward pointing unit
normal vectors, and $\bfa{e}_j$ is the unit vector in the $j$th direction.

The tensor $\tilde{\varkappa}(s)$ is bounded, symmetric, positive definite, as shown by \eqref{epk}, \eqref{k8} and \eqref{green} in the following manner \cite{homo02ye}.
First, let $w_i = \psi_i - x_i\,,$ for $i=1,\cdots, d\,.$ 
We have
\[\nabla w_i = \nabla \psi_i - \bfa{e}_i\,.\]
Second, \eqref{green} and \eqref{k8} leads to \begin{align}\label{psd}
\begin{split}
\bfa{e}_i \cdot \tilde{\varkappa}(s) \bfa{e}_j
&= \bfa{e}_i \cdot \tilde{\varkappa}(s) \langle  \nabla \psi_j  \rangle_K
 = \bfa{e}_i \cdot \langle \varkappa(\bfa{x},s) \nabla \psi_j  \rangle_K \\
 &=  \langle \bfa{e}_i \cdot \varkappa(\bfa{x},s) \nabla \psi_j  \rangle_K = 
 \langle (\nabla \psi_i - \nabla w_i) \cdot \varkappa(\bfa{x},s) \nabla \psi_j  \rangle_K\,.
\end{split}
\end{align}
Next, by \eqref{epk}, Divergence Theorem, the production rule and \eqref{epk0}, it holds that
\begin{align}\label{dvp}
\begin{split}
0 &= \int_{\partial K} w_i \bfa{n} \cdot \varkappa(\bfa{x},s) \nabla \psi_j \, \tu{d} \sigma \\	&= \int_{K} \div(w_i \varkappa(\bfa{x},s) \nabla \psi_j) \dx\\
& = \int_{K} w_i \div(\varkappa(\bfa{x},s) \nabla \psi_j) \, \dx + \int_{K} \nabla w_i \cdot \varkappa(\bfa{x},s) \nabla \psi_j \, \dx \\
& = \int_{K} \nabla w_i \cdot \varkappa(\bfa{x},s) \nabla \psi_j \, \dx\,.
\end{split}
\end{align}
Hence, \eqref{psd} becomes
\begin{align}\label{psd1}
\bfa{e}_i \cdot \tilde{\varkappa}(s) \bfa{e}_j
= \langle \nabla \psi_i \cdot \varkappa(\bfa{x},s) \nabla \psi_j  \rangle_K\,,
\end{align}
which also offers a different approach to calculate $\tilde{\varkappa}\,.$

\bigskip

To illustrate the effective coefficient calculation, we consider the Haverkamp model \cite{haverkamp, Spiridonov2019} possessing the following constitutive relation for the unsaturated hydraulic conductivity $\varkappa$:
\begin{equation}\label{steady}
	\varkappa(\bfa{x},p(\bfa{x})) = K_s(\bfa{x})\, K_r(p) = \kappa(\bfa{x}) \, \frac{1}{1+|p|}\,.
\end{equation}
Here, $\kappa(\bfa{x})$ denotes the permeability of soils, $K_r(p)$ represents the relative hydraulic conductivity, $K_s(\bfa{x})$ stands for the saturated hydraulic conductivity (which is mostly different from the permeability $\kappa(\bfa{x})$, see \cite{satper}).
All the multiscale heterogeneity is incorporated in $\kappa(\bfa{x})$ without regard to $p\,,$ and $\dfrac{1}{1+|p|}$ includes all the nonlinearity.  With this Haverkamp model \eqref{steady}, over each coarse cell $K \in \mathcal{T}^H\,,$ the effective permeability tensor is denoted by $\tilde{\kappa}$ and obtained from \eqref{k8}--\eqref{k8e}, where
\begin{equation}\label{efper}
\tilde{\varkappa}(\bfa{x},p(\bfa{x})) =  \tilde{\kappa}(\bfa{x}) \, \frac{1}{1+|p|}\,.
\end{equation}

In our two-dimensional case, on each coarse cell $K \in \mathcal{T}^H\,,$ the constant effective permeability tensor $\tilde{\kappa}$ from \eqref{steady} is of the following form
\begin{align}\label{k2d}
	\tilde{\kappa} = \left (
	\begin{matrix}
		\tilde{\kappa}_{11} & \tilde{\kappa}_{12}\\
		\tilde{\kappa}_{12} & \tilde{\kappa}_{22}
	\end{matrix}
	\right )\,.	
\end{align}

\subsubsection{Solving the homogenized Richards equation}\label{rep}

Finally, referring to the original equation \eqref{eq:original0}, we solve the following homogenized Richards equation on coarse grid with the effective coefficient $\tilde{\varkappa}$ \eqref{k8e}:
\beq
\label{eq:original0u}
\frac{\partial p(t,\bfa{x})}{\partial t} - \div (\tilde{\varkappa}(\bfa{x},p(t,\bfa{x}))\nabla p(t,\bfa{x}))
	= f(t,\bfa{x}) \  \textrm{in} \ (0,T] \times \Omega\,,
	\eeq
It is clear that this homogenized Richards equation is principally the same as the original one \eqref{eq:original0}, with the exception that scalars become symmetric and positive definite tensors.

For numerical simulations, we consider the time-dependent Richards equation \eqref{eq:original0} employing the Haverkamp model \eqref{steady}:
\beq
\label{eq:original00t}
\frac{\partial p(t,\bfa{x})}{\partial t} - \div\left ( \kappa(\bfa{x}) \, \frac{1}{1+|p|}\, \nabla p(t,\bfa{x}) \right )
= f(t,\bfa{x}) \  \textrm{in} \ (0,1] \times \Omega\,,
\eeq
where $\bfa{x}=(x_1,x_2)\,,$ $f(t,\bfa{x}) = \cos(\pi x_1)\sin(\pi x_2)\,,$ taking the initial condition $p(0,\bfa{x})= p_{0}=0$ in $\Omega$ and the Dirichlet boundary condition $p(t,\bfa{x})=0$ on $(0,T] \times \partial \Omega\,.$ 

After homogenization (addressed in Subsubsection \ref{ep}), we solve the homogenized Richards equation derived from \eqref{eq:original00t} on coarse grid with the effective permeability $\tilde{\kappa}$ \eqref{efper} as follows:
\beq
\label{rhomot}
\frac{\partial p(t,\bfa{x})}{\partial t} - \div\left (\tilde{\kappa}(\bfa{x}) \, \frac{1}{1+|p|}\, \nabla p(t,\bfa{x}) \right )
= f(t,\bfa{x}) \  \textrm{in} \ (0,1] \times \Omega\,,
\eeq
having the initial condition $p(0,\bfa{x})= p_{0}$ in $\Omega$ and the Dirichlet boundary condition $p(t,\bfa{x})=0$ on $(0,T] \times \partial \Omega\,.$ 

We denote
\[\tilde{a}(p,v;u) = \int_{\Omega} \tilde{\varkappa}(\bfa{x},u) \nabla p \cdot \nabla v \, \dx \,,\]
with a given $u \in V\,,$ for any $p,v \in V\,,$ and $\tilde{\varkappa}$ is from \eqref{k8e} as well as \eqref{efper}.  Let $P_H$ stand for the $L^2(\Omega)$ projection operator onto $V^H\,,$ and $p_0$ comes from \eqref{ini}.

The completely Picard discrete iterative scheme for the homogenized Richards equation \eqref{rhomot} reads in $\mathcal{T}^H$: starting with an initial $p_{c,0} = P_H p_0$ in $V^H\,,$ we guess $p^0_{c,s+1}\in V^H$ at the time step $(s+1)$th, 
and carry out iteration in $V^H$ to find $p^{n+1}_{c,s+1}\in V^H\,:$
\begin{align}\label{r1elhomoc}
	\left(\frac{p^{n+1}_{c,s+1}}{\tau} , v \right) + \tilde{a}(p^{n+1}_{c,s+1},v;p^n_{c,s+1}) 
	= \left(\frac{p_{c,s}}{\tau} , v \right) + (f_{s+1},v)\,,
\end{align}
for $n=0,1,2, \cdots \,,$ until meeting \eqref{pt} at some $\alpha$th Picard iteration, with all $v \in V^H\,.$
We use \eqref{pdata} to set the preceeding time data $p_{c,s+1} = p_{c,s+1}^{\alpha}$ to shift to the next temporal step in \eqref{r1ed}.
\bigskip

The finite element basis function at node $j$ (that is $\bfa{x}_j$) of the coarse grid $\mathcal{T}^H$ is denoted by $\varphi_j \in V^H\,.$  
The subscript $c$ can be removed from the following for the sake of simplicity.
At the time step $(s+1)$th, the approximate solution $p^{n+1}_{c,s+1}$ is
\begin{equation}\label{papprox}
	p^{n+1}_{c,s+1} \approx \sum_{j=1}^{N_v} p^{n+1}_{c,s+1}(\bfa{x}_j) \, \varphi_j = \sum_{j=1}^{N_v} p^{n+1}_{c,s+1,j} \, \varphi_j\,,
\end{equation}
where $p^{n+1}_{c,s+1,j}$ is the nodal point value of $p^{n+1}_{c,s+1}$ at the node $j\,,$ and $N_v$ is the total number of nodes on the coarse grid $\mathcal{T}^H\,.$  Let \begin{equation}\label{components}
\bfa{p}^{n+1}_{c,s+1} = (p^{n+1}_{c,s+1,1}, p^{n+1}_{c,s+1,2}, \cdots, p^{n+1}_{c,s+1,N_v})\,.
\end{equation}
Using the above form \eqref{papprox} and \cite{modifiedPicard}, we can write
\begin{align}\label{comkth}
	\begin{split}
		\varkappa(\bfa{x},p^n_{c,s+1}) &\approx \sum_{j=1}^{N_v} \varkappa(\bfa{x}_j,p^n_{c,s+1,j})\varphi_j = \sum_{j=1}^{N_v} \varkappa_j \, \varphi_j\,.
	\end{split}
\end{align}


With $N_v$ defined in \eqref{papprox} and for $i,j = \overline{1,N_v}\,,$ we let
\begin{align}\label{comstiffh}
\bfa{C}_{s+1} = \{C_{s+1,ij}\} \,, \quad 
\tilde{\bfa{A}}^n_{s+1} = \{\tilde{A}^n_{s+1,ij}\}\,, \quad
\bfa{b}_{s+1} = \{b_{s+1,i}\}\,, 
\end{align}
where
\begin{align}\label{comstiffh2}
\begin{split}
C_{s+1,ij} &= \int_{\Omega}  \varphi_j \cdot \varphi_i \, \dx \,,\\
\tilde{A}^n_{s+1,ij} &= \int_{\Omega} \tilde{\varkappa}(\bfa{x},p^n_{c,s+1}) \nabla \varphi_j \cdot \nabla \varphi_i \, \dx \,,\\
b_{s+1,i} &=   \int_{\Omega} \left(\frac{p_{s}}{\tau}  + f_{s+1}\right) \, \varphi_i \, \dx \,.
\end{split} 
\end{align}
Benefiting from \cite{rmrep1}, we write \eqref{r1elhomoc} in the following matrix form:
\begin{equation}\label{matrixf}
\bfa{C}_{s+1} \frac{\bfa{p}^{n+1}_{c,s+1}}{\tau} + \tilde{\bfa{A}}^n_{s+1} \bfa{p}^{n+1}_{c,s+1} =  \bfa{b}_{s+1}\,,
\end{equation}
where $\tilde{\bfa{A}}^n_{s+1}$ denotes the coarse-scale homogenized stiffness matrix, $\bfa{b}_{s+1}$ represents the right-hand side load vector, each of $\bfa{p}_{c,s+1}^n$ and $\bfa{p}^{n+1}_{c,s+1}$ is a vector of $N_v$ components defined in \eqref{components}.

\bigskip

\section{Deep learning for numerical homogenization}\label{dlhomo}


Regarding applications of the Richards equation \eqref{eq:original00t} as an unsaturated flow, on the one hand, the numerical homogenization gives us a quick solver to calculate the homogenized solutions efficiently and precisely.  On the other hand, there exist uncertainties in some local portions of the heterogeneous permeability field $\varkappa(\bfa{x},p(\bfa{x}))\,.$  To estimate the flow solution's uncertainties, several thousand forward simulations are required.  
We want to design a strategy that uses a given vast volume of simulation data and reduces direct calculating endeavor afterwards. For example, convolutional neural network (CNN) has recently been applied to numerical homogenization \cite{cnn-homo-elasticity19, cnn-homo-poro20}.  In this paper, we assume that the Picard iterative procedure \eqref{r1elhomoc} 
for the homogenized Richards equation \eqref{rhomot} stops at the $(s+1)$th time step and $(n+1)$th Picard step in \eqref{r1elhomoc}.  Then, deep neural networks (DNNs) are utilized to describe the link between stochastic permeability field $\kappa(\bfa{x})$ 
%
and expected numerical homogenization's macroscopic parameters, that is, coarse-scale effective permeability $\tilde{\kappa}(\kappa)$ \eqref{k2d}, homogenized stiffness matrix $\tilde{\bfa{A}}^{n}(\kappa)$ \eqref{comstiffh} (by data obtained from the $n$th Picard step) as well as the right-hand side vector $\bfa{b}(\kappa)$ \eqref{comstiffh2} of \eqref{r1elhomoc}.  Note that our stochastic permeability fields are generated using Karhunen-Lo\`eve expansion (KLE) \cite{cnn-mcmc20}. 
Subsequent to the above relation's establishment, for any random permeability field realization, we can feed it to the network 
and get relevant DNN-prediction of numerical homogenization's macroscopic parameters (to be specified in Subsection \ref{nnarch}), which will be applied to 
Subsubsection \ref{rep} for computing the predicted global coarse-scale solution of \eqref{rhomot}. 
For simplicity in the rest of this work, we will drop the subscript ($s+1$) where Picard steps are involved.


More specifically, in order to create the neural network, we are intrigued by building the nonlinear maps $g_{\kappa}\,,$ $g_{A}^{n}\,,$ and $g_{b}$ that respectively map the permeability field $\kappa$ to the desired macroscopic parameters (effective permeability $\tilde{\kappa}(\kappa)\,,$ coarse-scale homogenized stiffness matrix $\tilde{\bfa{A}}^{n}(\kappa)\,,$ and right-hand side vector $\bfa{b}(\kappa)$ of \eqref{r1elhomoc}) as follows:
\begin{align}
g_{\kappa}: \kappa &\to \tilde{\kappa}(\kappa)\,,\label{gk}\\
g_{A}^{n}: \kappa &\to \tilde{\bfa{A}}^{n}(\kappa)\,,\label{gm}\\
g_{b}: \kappa &\to \bfa{b}(\kappa)\,.\label{gb}
\end{align}	

This study's purpose is to take advantage of deep learning to create rapid approximations of the above expected macroscopic parameters related to the permeability field $\kappa$'s uncertainties, allowing us to solve the homogenized Richards equation \eqref{rhomot} quickly and accurately.  For each permeability field $\kappa\,,$ its images can be calculated thanks to the maps $g_{\kappa}\,,$ $g_{A}^{n}\,,$ and $g_{b}\,.$ 
We utilize these forward computations as training samples for constructing deep neural networks, which will approximate the associated maps, that is, 
\begin{align}\label{nnbm}
\begin{split}
\mathcal{N}_{\kappa}  (\kappa) &\approx g_{\kappa} (\kappa)\,,\\
\mathcal{N}_{A}^{n}  (\kappa) &\approx g_{A}^{n} (\kappa)\,,\\
\mathcal{N}_{b}  (\kappa) &\approx g_{b} (\kappa)\,.
\end{split}
\end{align}

Within such neural networks, the input is permeability field $\kappa\,;$ whereas the predicted outputs are macroscopic parameters, expressly, effective permeability $\tilde{\kappa}(\kappa)\,,$ coarse-scale homogenized stiffness matrix $\tilde{\bfa{A}}^{n}(\kappa)\,,$ and right-hand side vector $\bfa{b}(\kappa)$ of \eqref{r1elhomoc}. Upon constructing these neural networks, they can be used for any new permeability realization $\kappa\,,$ to predict output macroscopic parameters $\tilde{\kappa}_{\tu{pred}}\,,$ $\tilde{\bfa{A}}^{n}_{\tu{pred}}\,,$ and $\bfa{b}_{\tu{pred}}\,,$ which will be defined below.


\subsection{Network Architecture}\label{nnarch}

The neurons in most deep neural networks can be arranged into three categories: the input layer, then one or several hidden layers, and finally the output layer.  In the majority of deep neural network architectures, these layers are organized in a chain pattern, where each layer acts as a function of the layer before it \cite{rpdnn, dnngood}.
Thus, at least one hidden layer exists in every neural network (also called artificial neural network); otherwise, a neural network is not what it is. Deep neural networks are described as having at least two or usually many more hidden layers.

Normally, an $L$-layer feedforward and fully connected deep neural network $\mathcal{N}^L$ can be represented in the following form \cite{dnngms19,rpdnn, dnngood}:
\begin{equation}\label{lnn}
\mathcal{N}(z;\theta) = \mathcal{N}^L(z;\theta) = \sigma^L(\bfa{W}^{L} \sigma^{L-1}(\cdots \sigma^2(\bfa{W}^2 \sigma^1(\bfa{W}^1 z + \bfa{c}^1) + \bfa{c}^2) \cdots) + \bfa{c}^{L})\,,
\end{equation}
where $z$ is the input layer, $L$ is the network's depth as the number of layers (excluding the input layer), $\theta:=(\bfa{W}^1,\bfa{W}^2, \cdots, \bfa{W}^{L}, \bfa{c}^1, \bfa{c}^2, \cdots, \bfa{c}^{L})$ with the  bias vectors $\bfa{c}^l$ as well as the weight matrices $\bfa{W}^l$ (for $l=1,\cdots,L$), and $\sigma^l$ are the activation functions (which are often nonlinear as SELU or ReLU for the input layer or hidden layers).

More specifically, the input layer of this neural network is
\[\mathcal{N}^0(z) = z\,.\]
For $1 \leq l < L\,,$ the $l$th hidden layer of the network is
\[\mathcal{N}^l(z) = \sigma^l (\bfa{W}^{l}\, \mathcal{N}^{l-1}(z)  + \bfa{c}^{l})\,.\]
The output layer (without using any activation function) is
\[\mathcal{N}^L(z) = \bfa{W}^{L}\, \mathcal{N}^{L-1}(z)  + \bfa{c}^{L}\,.\]

\begin{figure}[ht!]
	\centering
	\includegraphics[width=0.6\linewidth]{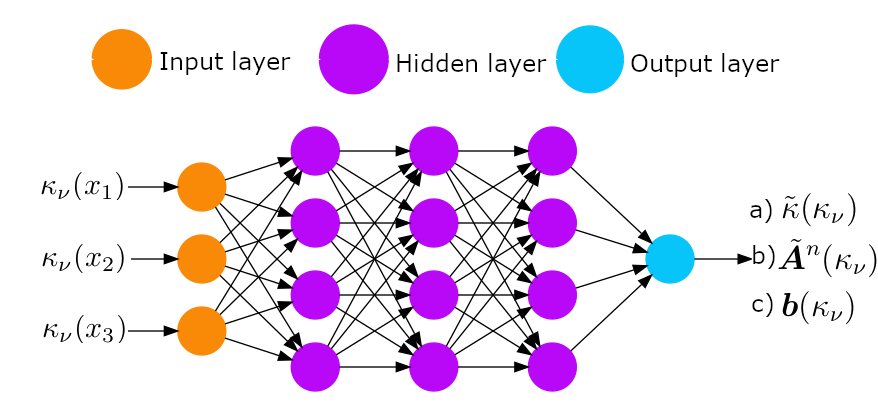}
	\caption{Illustration of a deep neural network.}
	\label{dnn_illustrate}
\end{figure}

A neural network of this type is utilized for approximating the expected output $y$ (also known as target,  desired output, correct output, or training value).  By solving the following optimization problem, we aim at determining $\theta^*\,:$ 
\begin{equation}\label{opt}
\theta^{*} = \argmin_{\theta} \mathcal{L}(\theta)\,.
\end{equation}
Here, $\mathcal{L}(\theta)$ stands for the loss function and represents the difference's measurement between the image (predicted output) of the network $\mathcal{N}(z;\theta)$'s input and the target $y$ in a training set consisting of training pairs $(z_{\nu},y_{\nu})\,,$ which are also called training samples, or training examples, or training patterns.  In this work, the loss function is the mean squared error:
\begin{equation}\label{loss}
\mathcal{L}(\theta) = \frac{1}{N} \sum_{\nu=1}^N \|y_{\nu} - \mathcal{N}(z_{\nu}; \theta)\|^2_2\,,
\end{equation}
where $N$ represents how many training samples there are.  Figure \ref{dnn_illustrate} depicts an example of a deep neural network (to be explained below).

%
 
Consider that we are given a collection of distinct permeability realizations $\{\kappa_{\nu}\}_{\nu=1}^N$ over the whole domain $\Omega\,.$  Within our neural network, the input $z_{\nu} = \kappa_{\nu} \in \mathbb{R}^m$ is a vector (named input vector or feature vector) representing the permeability field's values as features $\kappa_{\nu}(\bfa{x}_1), \kappa_{\nu}(\bfa{x}_2),\cdots,\kappa_{\nu}(\bfa{x}_m)$, where $\bfa{x}_1, \bfa{x}_2,\cdots, \bfa{x}_m \in \mathbb{R}^d$ are the spatial points in the domain $\Omega\,,$ with $m=256$ and $d=2$ in our simulations.  Recall that our stochastic permeability fields are produced by Karhunen-Lo\`eve expansion (KLE) \cite{cnn-mcmc20}, to be explained in Section \ref{numer} on numerical experiments.
The correct output $y_{\nu}$ 
is a vector
of a vectorized coarse-grid effective permeability tensor $\tilde{\kappa}({\kappa_{\nu}})\,,$ or homogenized stiffness matrix $\tilde{\bfa{A}}^{n}({\kappa}_{\nu})\,,$ or right-hand side vector $\bfa{b}({\kappa}_{\nu})$ of \eqref{r1elhomoc} by numerical homogenization in Subsection \ref{sec:method}.  
The deep neural networks $\mathcal{N}_{\kappa}(z,\theta_{\kappa})\,,$ $\mathcal{N}^n_{A}(z,\theta_{A})\,,$ and $\mathcal{N}_{b}(z,\theta_{b})$ described in \eqref{nnbm} and \eqref{lnn} are trained using the training pairs $(z_{\nu}, y_{\nu})$ by minimizing the loss function $\mathcal{L}$ with regard to $\theta$ as the network parameter, so that the trained neural networks $\mathcal{N}_{\kappa}(z,\theta^*_{\kappa})\,,$ $\mathcal{N}^n_{A}(z,\theta^*_{A})\,,$ and $\mathcal{N}_{b}(z,\theta^*_{b})$ can respectively estimate the functions $g_{\kappa}\,,$ $g_{A}^{n}\,,$ and $g_{b}$ in \eqref{gk}--\eqref{gb} on coarse grid.  After creating the neural networks, we utilize them for any new permeability field $\kappa_{N+\iota}$ ($\iota = 1,\cdots,M$) to quickly predict the coarse-scale outputs, which are effective permeability $\tilde{\kappa}_{\tu{pred}}$ by
\begin{equation}\label{kpred}
\tilde{\kappa}_{\tu{pred}}(\kappa_{N+\iota}) = \mathcal{N}_{\kappa}(\kappa_{N+\iota};\theta^*_{\kappa}) \approx g_{\kappa}(\kappa_{N+\iota}) =  \tilde{\kappa}(\kappa_{N+\iota})\,,
\end{equation}
and homogenized stiffness matrix $\tilde{\bfa{A}}^{n}_{\tu{pred}}$ by 
\begin{equation}\label{mpred}
\tilde{\bfa{A}}_{\tu{pred}}^{n}(\kappa_{N+\iota}) = \mathcal{N}_{A}^{n}(\kappa_{N+\iota};\theta^*_A) \approx g_{A}^{n}(\kappa_{N+\iota}) =  \tilde{\bfa{A}}^{n}(\kappa_{N+\iota})\,,
\end{equation}
and right-hand side vector $\bfa{b}_{\tu{pred}}$ of \eqref{r1elhomoc} by
\begin{equation}\label{bpred}
\bfa{b}_{\tu{pred}}(\kappa_{N+\iota}) = \mathcal{N}_{b}(\kappa_{N+\iota};\theta^*_b) \approx g_{b}(\kappa_{N+\iota}) =  \bfa{b}(\kappa_{N+\iota})\,.
\end{equation}




\subsection{Network-based numerical homogenization solver}\label{nnsolve}
Significantly, upon designing the neural networks, following the matrix form \eqref{matrixf},
at the final ($s+1$)th temporal step and $(n+1)$th Picard step in \eqref{r1elhomoc}, we can find the predicted global coarse-scale solution $p^{n+1}_{c,\tu{pred}}$ of the homogenized Richards equation \eqref{rhomot} by solving for the coarse-scale coefficient vector 
$\bfa{p}^{n+1}_{c,\tu{pred}}$ \eqref{components} from the linear algebraic system
\begin{equation}\label{apred}
\bfa{C}_{s+1} \frac{\bfa{p}^{n+1}_{c,\tu{pred}}}{\tau} + \tilde{\bfa{A}}^n_{\tu{pred}} \ \bfa{p}^{n+1}_{c,\tu{pred}} = \bfa{b}_{\tu{pred}}\,.
\end{equation}
Here, $\tilde{\kappa}_{\tu{pred}}\,,$ $\tilde{\bfa{A}}^{n}_{\tu{pred}}\,,$ and $\bfa{b}_{\tu{pred}}$ were predicted in \eqref{kpred}, \eqref{mpred}, and \eqref{bpred}, respectively.

\section{Numerical Examples}\label{numer}

\begin{figure}[H]
	\begin{center}
	  \includegraphics[width=\linewidth]{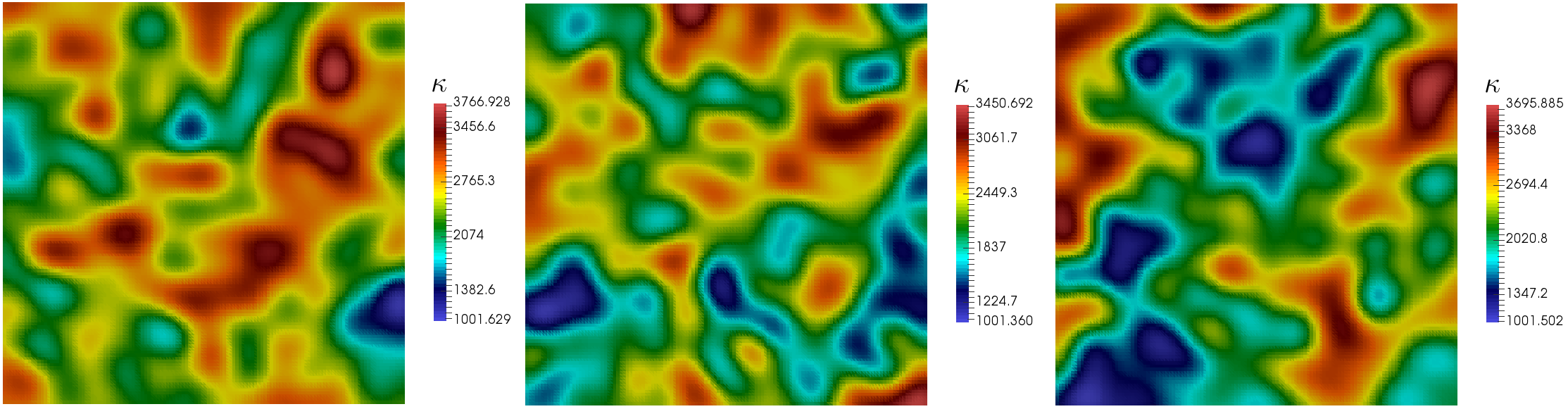}
	\end{center}
	\caption{Examples of three permeability fields $\kappa$ in range [1000, 4200].}
	\label{fig:kappa}
\end{figure}

Within this part, utilizing our developed deep learning strategy, several numerical results are shown for the predicted macroscopic parameters and thus solutions of the Richards equation's numerical homogenization. More specifically, we investigate the accuracy of deep neural networks (DNNs). Numerical experiments are carried out for the steady-state and time-dependent homogenized Richards equation \eqref{rhomot} derived from the Haverkamp model \eqref{eq:original00t} based on the original Richards equation \eqref{eq:original0}.  The considered permeability fields $\kappa(\bfa{x})$ in \eqref{eq:original00t} are heterogeneous, high-contrast, and stochastic (to be specified in Subsection \ref{sper}). Some examples of our generated permeabilities are presented in Figure  \ref{fig:kappa}.  On coarse grid, for each permeability field $\kappa(\bfa{x})\,,$ we investigate two DNN approaches: (1) prediction of effective permeability tensor 
and (2) prediction of homogenized stiffness matrix 
(for steady-state case) 
together with prediction of right-hand side vector 
of \eqref{r1elhomoc} (for time-dependent case).  

In simulations, our computational domain $\Omega = [0,1] \times [0,1]$ is heterogeneous and includes uncertainties. We use $8 \times 8$ square coarse grid as Fig.~\ref{fig:mesh} (of edge size $H=1/8$) and its refinement fine grid as Fig.~\ref{fig:fine1block} of $128 \times 128$ squares with two triangles (of size $h=\sqrt{2}/128$) per square fine cell.  The starting guess for pressure of the Picard iteration process is $0\,.$
Over all numerical experiments, the Picard iterative termination criterion is $\delta_0 = 10^{-6}\,,$ which ensures this linearization procedure's convergence, and we assume that the final Picard step is $(n+1)$th. The maximum number of Picard iterations is 4, for either the steady-state or time-dependent case.  



\subsection{Stochastic permeability}\label{sper}
    
In order to parameterize the permeability fields' uncertainties, we employ the Karhunen-Lo\`eve expansion (KLE) \cite{aarnes2008mixed,cnn-mcmc20,cnn-homo-poro20, ganis2008stochastic,kle2}.  Let  $\Upsilon (\bfa{x},\omega)$ be a random process.  Its covariance function is denoted by $R(\bfa{x},\hat{\bfa{x}})\,,$ with the position vectors $\bfa{x}=(x_1,x_2)\,, \ \hat{\bfa{x}}=(\hat{x}_1,\hat{x}_2)$ in $\Omega\,.$
We choose a special orthonormal basis $\{\phi_k(\bfa{x})\}$ in the Hilbert space $L^2(\Omega)$ such that
$\phi_k(\bfa{x})$ are the eigenfunctions of $R(\bfa{x},\hat{\bfa{x}})$ in the spectral problem
\begin{equation}\label{spec}
\int_{\Omega} R(\bfa{x},\hat{\bfa{x}})\phi_k(\hat{\bfa{x}})\tu{d} \hat{x} = \lambda_k \phi_k(\bfa{x})\,, \quad k = 1,2,\cdots\,,
\end{equation}
where $\lambda_k$ are the eigenvalues.  With these $\phi_k$ and $\lambda_k$ obtained from \eqref{spec}, $\Upsilon (\bfa{x},\omega)$ can be expressed in the following form called Karhunen-Lo\`eve expansion \cite{kle2}: 
\begin{equation}\label{randfield}
    \Upsilon (\bfa{x},\omega) = \sum^{\infty}_{k=1} \sqrt{\lambda _k} \zeta_k(\omega) \phi_k({\bfa{x})}\,,
\end{equation}
where $\zeta_k(\omega)$ are the random variables to be found.  

In experiments, regarding \eqref{randfield}, only $N_{\Upsilon}$ leading terms (determined by $\lambda_k$'s size) need to be kept to preserve most of $\Upsilon (\bfa{x},\omega)$'s energy, and thus, \[\Upsilon (\bfa{x},\omega) = \sum^{N_{\Upsilon}}_{k=1} \sqrt{\lambda _k} \zeta_k(\omega) \phi_k({\bfa{x})}\,.
\]
About the covariance function $R(\bfa{x},\hat{\bfa{x}})\,,$ it is supposed to be of the form
\begin{equation}\label{co-matrix}
    R(\bfa{x},\hat{\bfa{x}})=\sigma^2 \text{exp} \Bigg(-\sqrt{\frac{|x_1-\hat{x}_1|^2}{\eta^2_1}+\frac{|x_2-\hat{x}_2|^2}{\eta^2_2}}\Bigg),
\end{equation}    
abbreviated by $\kappa(\bfa{x})\,,$ where $\eta_1 = 0.2\,, \eta_2=0.2$ are correlation lengths in each spatial direction, and $\sigma^2=2$ is the variance.  Then, each of our stochastic permeability fields $\kappa(\bfa{x},\omega)$ can be represented by 
\begin{equation}\label{kleform}
    \kappa(\bfa{x},\omega) = \text{exp}(a_k (m(\bfa{x},\omega )),
\end{equation}
where $m(\bfa{x},\omega ) = m(\Upsilon (\bfa{x},\omega))$ is the heterogeneous porosity, and $a_k > 0\,.$  Recall that some examples of our generated permeabilities are presented in Figure  \ref{fig:kappa}.

Alternative ways to build a stochastic permeability field can be employed without affecting our described deep learning methodology.  The reason is that a separate neural network is created for each of the macroscopic parameters (effective permeability tensor, homogenized stiffness matrix, and right-hand side vector of \eqref{r1elhomoc}) as follows.



\subsection{Overview of deep learning for numerical homogenization}\label{numdnn}

We use the numerical homogenization technique described in Subsection \ref{sec:method} to produce the coarse-scale sample pairs (by data obtained from the $n$th Picard step of \eqref{r1elhomoc}).  The stochastic permeability field $z = \kappa$ is viewed as the neural network's input, whereas the desired outputs are coarse-scale effective permeability tensor $\tilde{\kappa}$ \eqref{k2d}, homogenized stiffness matrix $\tilde{\bfa{A}}^n$ \eqref{comstiffh} (for steady-state case) along with right-hand side vector $\bfa{b}$ \eqref{comstiffh2} from \eqref{r1elhomoc} (for time-dependent case).  We arbitrarily split up these sample pairs into the training set (also known as training data, training dataset, or learning set) and testing set.  
The training set is formed by sample pairs produced by a high number $N$ of random permeability realizations $\kappa_1, \kappa_2, \cdots, \kappa_N\,.$  On the other hand, the testing set consists of the sample pairs constructed by the remaining $M$ random permeability realizations $\kappa_{N+1},\kappa_{N+2}, \cdots, \kappa_{N+M}$, for assessing the trained networks's predictive ability.  We note that, in order to create a unique deep neural network, an optimization problem is addressed for each of the macroscopic parameters ($\tilde{\kappa}\,,$  $\tilde{\bfa{A}}^n\,,$ and $\bfa{b}$) 
by minimizing the loss function established by the training set.  The network architectures used to train such macroscopic parameters 
are summarized as follows (with $m = 256 = 16 \times 16$ and $\bfa{x}_1,\bfa{x}_2,\cdots, \bfa{x}_m \in \mathbb{R}^2$):
\begin{itemize}
\item[a)] For the effective permeability $\tilde{\kappa}\,,$ we generate a network $\mathcal{N}_{\kappa}$ using	
\begin{itemize} 
\item Input: vectorized permeability field $\kappa_{\nu}$ 
of values $\kappa_{\nu}(\bfa{x}_1)\,, \kappa_{\nu}(\bfa{x}_2)\,, \cdots, \kappa_{\nu}(\bfa{x}_m)\,,$ 
\item Desired output: vectorized effective permeability $\tilde{\kappa}(\kappa_{\nu})$ of 256 entries, 
\item Loss function \eqref{loss}: mean squared error $\dd  \frac{1}{N} \sum_{\nu=1}^N \|\tilde{\kappa}(\kappa_{\nu}) - \mathcal{N}_{\kappa}(\kappa_{\nu}; \theta_{\kappa})\|^2_{\mathbb{L}^2(\Omega)}\,,$
\item Activation function: SELU (Scaled Exponential Linear Unit) activation function for the first input layer, then ReLU (Rectified Linear Unit) activation function for all hidden layers, no activation function at the last output layer,
\item DNN structure: 3 hidden layers, each layer comprises 128-356 neurons,
\item Kernel initializer: normal, \item Training optimizer: Adam.
\end{itemize}
\end{itemize}

\bigskip

\begin{itemize}
\item[b)] For the coarse-grid homogenized stiffness matrix  $\tilde{\bfa{A}}^{n}\,,$ we create a network $\mathcal{N}^n_{A}$ employing	
\begin{itemize} 
\item Input: vectorized permeability field $\kappa_{\nu}$ 
of values $\kappa_{\nu}(\bfa{x}_1)\,, \kappa_{\nu}(\bfa{x}_2)\,, \cdots, \kappa_{\nu}(\bfa{x}_m)\,,$ 
\item Desired output: vectorized coarse-scale homogenized stiffness matrix $\tilde{\bfa{A}}^{n}(\kappa_{\nu})$ of 375 entries,  
\item Loss function \eqref{loss}: mean squared error $\dd  \frac{1}{N} \sum_{\nu=1}^N \|\tilde{\bfa{A}}^{n}(\kappa_{\nu}) - \mathcal{N}_{A}^{n}(\kappa_{\nu}; \theta_A)\|^2_{\mathbb{L}^2(\Omega)}\,,$
\item Activation function: SELU activation function for the first input layer, then ReLU activation function for all hidden layers, no activation function at the last output layer,
\item DNN structure: 3 hidden layers, each layer consists of 384-512 neurons,
\item Kernel initializer: normal, \item Training optimizer: Adam.
\end{itemize}
\end{itemize}

\begin{itemize}
\item[c)] For the right-hand side vector $\bfa{b}$ of \eqref{r1elhomoc}, we generate a network $\mathcal{N}_{M}^{n}$ having	
\begin{itemize} 
\item Input: vectorized permeability field $\kappa_{\nu}$ 
of values $\kappa_{\nu}(\bfa{x}_1)\,, \kappa_{\nu}(\bfa{x}_2)\,, \cdots, \kappa_{\nu}(\bfa{x}_m)\,,$ 
\item Desired output: right-hand side vector $\bfa{b}(\kappa_{\nu})$ of \eqref{r1elhomoc} with 81 entries, 
\item Loss function \eqref{loss}: mean squared error $\dd  \frac{1}{N} \sum_{\nu=1}^N \|\bfa{b}(\kappa_{\nu}) - \mathcal{N}_{b}(\kappa_{\nu}; \theta_b)\|^2_{\bfa{L}^2(\Omega)}\,,$
\item Activation function: SELU activation function for the first input layer, then ReLU activation function for all hidden layers, no activation function for the last output layer,
\item DNN structure: 3 hidden layers, each layer contains 256-384 neurons,
\item Kernel initializer: normal, 
\item Training optimizer: Adam.
\end{itemize}
\end{itemize}


We utilize the activation functions SELU (Scaled Exponential Linear Unit) \cite{selu} and ReLU (Rectified Linear Unit) \cite{ReLU28} because of their proved enablement in training deep neural networks without any vanishing gradient problem. 
Moreover, ReLU has the simplest derivative out of all nonlinear activation functions.  SELU is less common than ReLU, but SELU is very promising thanks to the fact that as opposed to ReLU, SELU cannot die.  By definition, SELU induces self-normalization neural network (SNN), and neuronal activations in the SNN automatically reach a zero mean and unit standard deviation.
SELU alone learns more quickly and effectively than other activation functions.  
%
The optimizer Adam (Adaptive Moment Estimation) is a further extension of stochastic gradient descent (SGD) 
and is widely employed in neural network training \cite{ada29}.  
%
%
Over all experiments, we use Python API Tensorflow and Keras \cite{keras30} to train our neural network.

Upon training, a neural network is constructed, and we can apply it to output's predicting provided a new
input. For the network to be helpful, the predictions must be accurate.  Through all experiments, to study the predictive ability of our established network, we employ $M$ sample pairs, which were not touched in the network training. With respect to these sample pairs as the testing set, we compare the predicted output with the target (correct output) and calculate their difference using some appropriate metric (to assess how well the model is working \cite{keras30}).

Regarding error formulas, over the training dataset on coarse grid, the error computation benefits from the relative $l^2$ error, namely, root mean square error (RMSE) \cite{nlnlmc31,rl2e1, rl2e2, rl2e3}:
%
\begin{equation}\label{rl2e}
RMSE = \sqrt{\frac{\sum_{\nu =1}^{N}| \hat{y}_{\nu} - y_{\nu}|^2}{\sum_{\nu = 1}^N |y_{\nu}|^2}}\,,
\end{equation}
where $N$ is the number of training pairs, 
$y_{\nu}$ is the target, and $\hat{y}_{\nu}$ is the predicted value for the sample pair $(z_{\nu}, y_{\nu})$ defined right before \eqref{loss}.  
%
%
%

On $M$ testing samples, the relative $l^2$ errors between the targets ($\tilde{\kappa}\,,\tilde{\bfa{A}}^{n}\,, \bfa{b}$) and their corresponding predicted coarse-scale effective permeability $\tilde{\kappa}_{\tu{pred}}$ \eqref{kpred}, homogenized stiffness matrix $\tilde{\bfa{A}}^{n}_{\tu{pred}}$ \eqref{mpred}, and right-hand side vector $\bfa{b}_{\tu{pred}}$ \eqref{bpred} of \eqref{r1elhomoc} are respectively
\begin{equation}\label{Cer}
\begin{split}
e^{\tilde{\kappa}}_{l^2}(\kappa_{N+ \iota})&=\frac{|\tilde{\kappa}_{\tu{pred}} (\kappa_{N+ \iota})-\tilde{\kappa}(\kappa_{N+ \iota})|}{|\tilde{\kappa}(\kappa_{N+ \iota})|}\,, \\
e^{\tilde{A}}_{l^2}(\kappa_{N+ \iota})&=\frac{|\tilde{\bfa{A}}^{n}_{\tu{pred}}(\kappa_{N+ \iota})-\tilde{\bfa{A}}^{n}(\kappa_{N+ \iota})|}{|\tilde{\bfa{A}}^{n}(\kappa_{N+ \iota})|}\,,\\
e^{b}_{l^2}(\kappa_{N+ \iota})&=\frac{|\bfa{b}_{\tu{pred}}(\kappa_{N+ \iota})-\bfa{b}(\kappa_{N+ \iota})|}{|\bfa{b}(\kappa_{N+ \iota})|}\,,
\end{split}
\end{equation}
where $\iota = 1,\cdots,M\,.$ Here, the norms were defined in Subsection \ref{sec:model}, the first two formulas has the Frobenius norm (or entrywise 2-norm) of matrix, and the last formula possesses the Euclidean norm (or 2-norm) of vector.
The formulas in \eqref{Cer} are also called metrics. 
%

More importantly, to evaluate the accuracy of our proposed DNN approaches for solutions, we employ the following relative $L^2$ and $H^1$ (energy) errors:
\begin{equation}\label{l2a}
\begin{split}
e^E_{L^2}=\frac{||p^E_{c,\tu{pred}}(\kappa_{N+\iota})-p_c(\kappa_{N+\iota})||_{L^2(\Omega)}}{||p_c(\kappa_{N+\iota})||_{L^2(\Omega)}},\quad
e^E_{H^1}=\frac{||\nabla p^E_{c,\tu{pred}}(\kappa_{N+\iota})-\nabla p_c(\kappa_{N+\iota})||_{\bfa{L}^2(\Omega)}}{||\nabla p_c(\kappa_{N+\iota})||_{\bfa{L}^2(\Omega)}}\,, \\
e^A_{L^2}=\frac{||p^A_{c,\tu{pred}}(\kappa_{N+\iota})-p_c(\kappa_{N+\iota})||_{L^2(\Omega)}}{||p_c(\kappa_{N+\iota})||_{L^2(\Omega)}},\quad
e^A_{H^1}=\frac{||\nabla p^A_{c,\tu{pred}}(\kappa_{N+\iota})-\nabla p_c(\kappa_{N+\iota})||_{\bfa{L}^2(\Omega)}}{||\nabla p_c(\kappa_{N+\iota})||_{\bfa{L}^2(\Omega)}}\,.
\end{split}
\end{equation}
In these expressions, $p_{c}$ is the homogenized solution of \eqref{matrixf}, $p^E_{c,\tu{pred}}$ is the solution for \eqref{apred} obtained by DNN-prediction of effective permeability, $p^A_{c,\tu{pred}}$ is the solution to \eqref{apred} derived by DNN-prediction of homogenized stiffness matrix (for the steady-state problem) or homogenized stiffness matrix together with right-hand side vector of \eqref{r1elhomoc} (for the time-dependent problem).  All these solutions are achieved from the last Picard step $(n+1)$th (and at the last time step $(s+1)$th in the time-dependent case), where the starting guess is $0$ for the pressure of Picard iterative procedure.

In experiments, we generate 6000 realizations of random heterogeneous permeability fields.
Recall that some samples of permeability fields are presented in Figure~\ref{fig:kappa}. These 6000 permeability realizations are then employed to produce 6000 sample pairs. 
As a typical practice, among these 6000 sample pairs, we arbitrarily pick $N=5000$ sample pairs as training samples (where 20\% of them are for validation) and use the rest $M=1000$ sample pairs as testing samples.

We should point out that raw data are not often fed directly to neural networks nor statistical models.
Usually, data need to be prepared in order to ease the network optimization procedure and increase the chances of getting favorable results. In particular, before training a neural network model, data processing like standardization or normalization are used to rescale input and output variables.
%
When a dataset is standardized, the value distribution is scaled to have a mean of 0 and a standard deviation of 1. Whereas, normalization is the procedure of forcing data to fit within a given range, which is $[-1,1]$ or $[0,1]\,,$ and is generally determined by the activation function being used.     

In this paper, only data normalization is considered.  Specifically, we are interested in a fairly simple normalizing technique, called reciprocal normalization, 
%
which normalizes numbers to the range -1 to 1, having the respective normalization and denormalization formulas
\begin{equation}
   x'= 2 \frac{x  - \tu{min}}{\tu{max} - \tu{min}} -1\,, \qquad       x = \bigg[\frac{x'+1}{2}(\text{max}-\text{min})\bigg]+\text{min}\,,
\end{equation}
Here, $x'$ is the normalized value, $x$ is the original (or denormalized) value, $\text{min}$ and $\text{max}$ are respectively the minimum and maximum observable values of the given data.  Note that such normalization method is applied to all input and target variables in this study. In simulations, our dataset's normalizing procedure benefits from the scikit-learn object MinMaxScaler. 

\subsection{Experiment for steady-state Richards equation}\label{ssc}

We first consider in this section the steady-state Richards equation derived from \eqref{eq:original00t} as a special case of \eqref{eq:original0}: find $p \in V$ such that
\beq
\label{eq:original00}
- \div \left(\kappa(\bfa{x}) \, \frac{1}{1+|p|}  \nabla p \right) = f \  \textrm{in } \Omega\,,
	\eeq
with the Dirichlet boundary condition $p(\bfa{x})=0$ on $\partial \Omega\,$ and $f=1\,.$  



\begin{figure}[H]
\centering
\begin{subfigure}{0.31\textwidth}
\centering
\footnotesize{Epochs -- 150 / Batch size -- 4}\\
\includegraphics[width=0.99\linewidth]{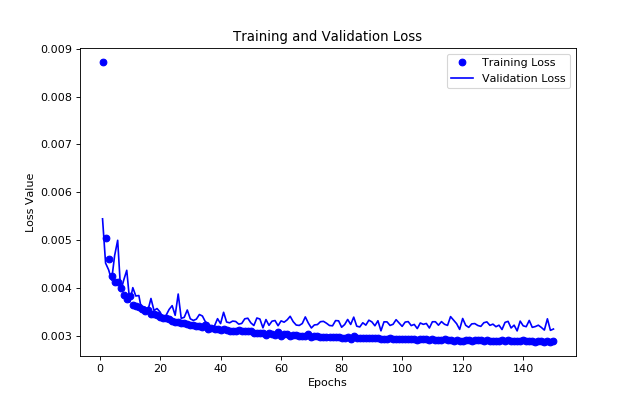}
\includegraphics[width=0.99\linewidth]{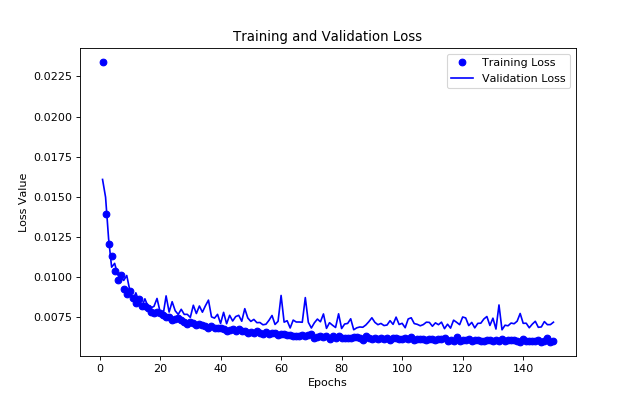}
\end{subfigure}
\begin{subfigure}{0.31\textwidth}
\centering
\footnotesize{Epochs -- 150 / Batch size -- 16}\\
\includegraphics[width=0.99\linewidth]{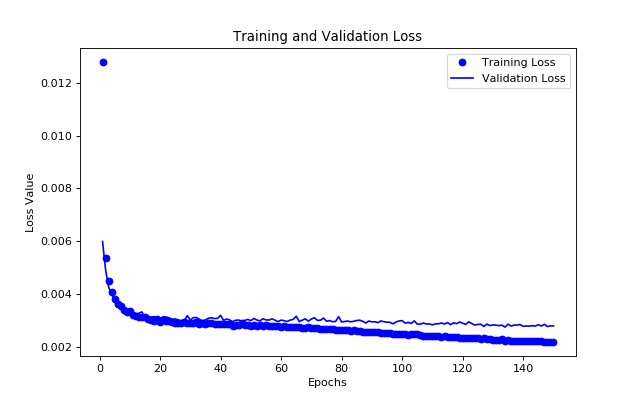}
\includegraphics[width=0.99\linewidth]{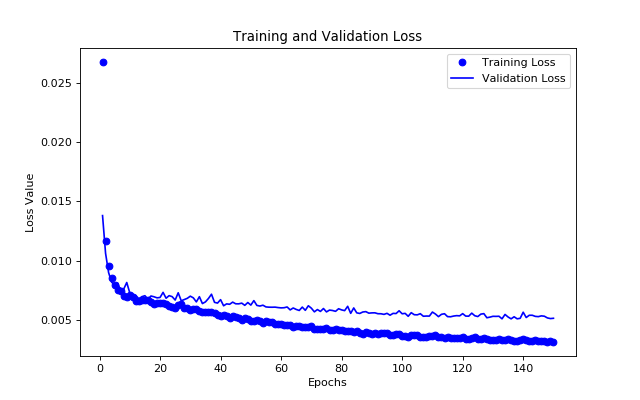}
\end{subfigure}
\begin{subfigure}{0.31\textwidth}
\centering
\footnotesize{Epochs -- 150 / Batch size -- 64}\\
\includegraphics[width=0.99\linewidth]{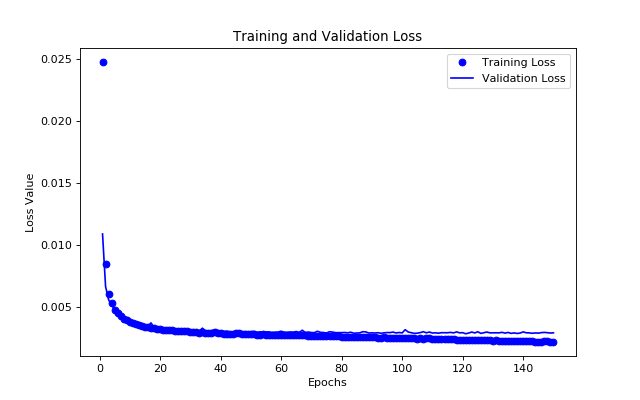}
\includegraphics[width=0.99\linewidth]{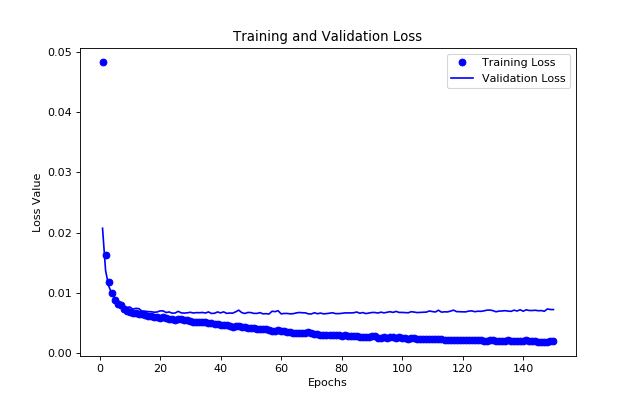}
\end{subfigure}
\begin{subfigure}{0.31\textwidth}
\centering
\footnotesize{Epochs -- 300 / Batch size -- 4}\\
\includegraphics[width=0.99\linewidth]{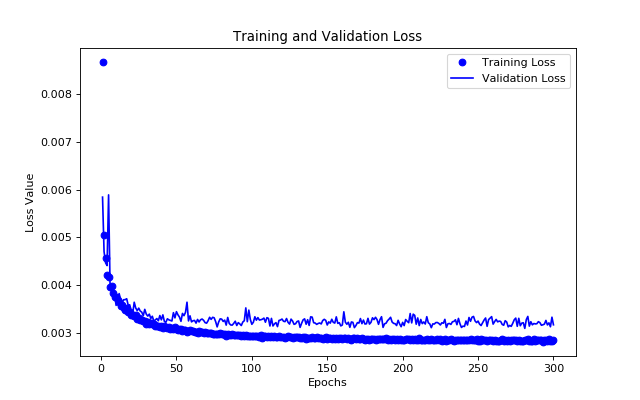}
\includegraphics[width=0.99\linewidth]{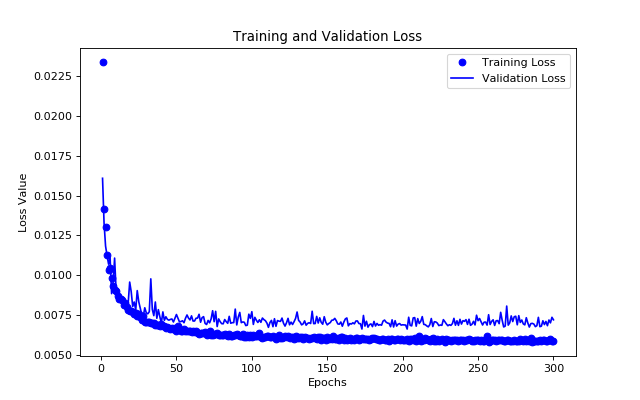}
\end{subfigure}
\begin{subfigure}{0.31\textwidth}
\centering
\footnotesize{Epochs -- 300 / Batch size - 16}\\
\includegraphics[width=0.99\linewidth]{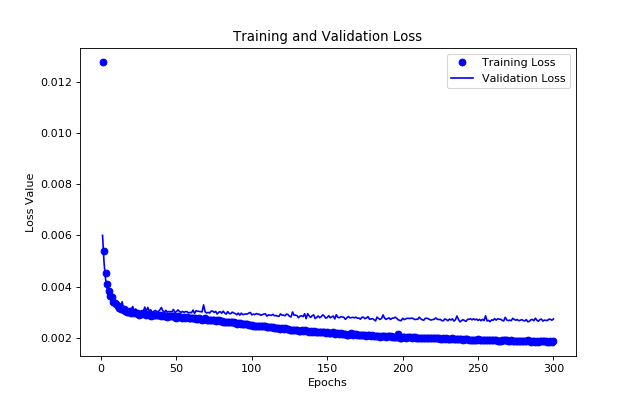}
\includegraphics[width=0.99\linewidth]{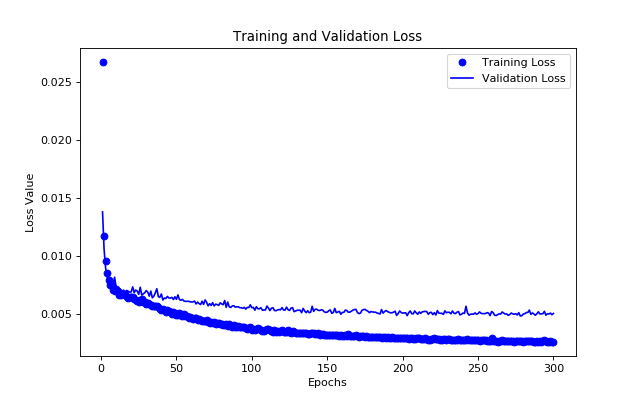}
\end{subfigure}
\begin{subfigure}{0.31\textwidth}
\centering
\footnotesize{Epochs -- 300 / Batch size - 64}\\
\includegraphics[width=0.99\linewidth]{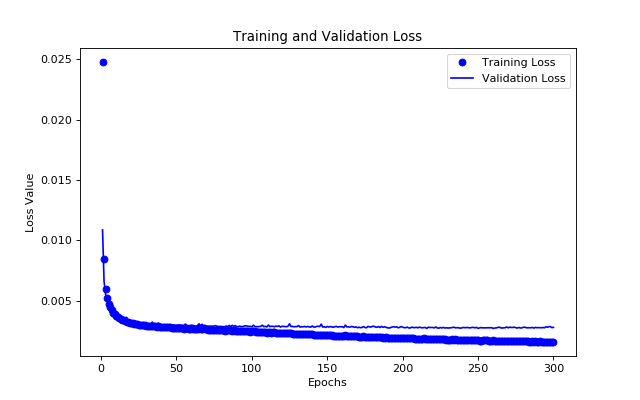}
\includegraphics[width=0.99\linewidth]{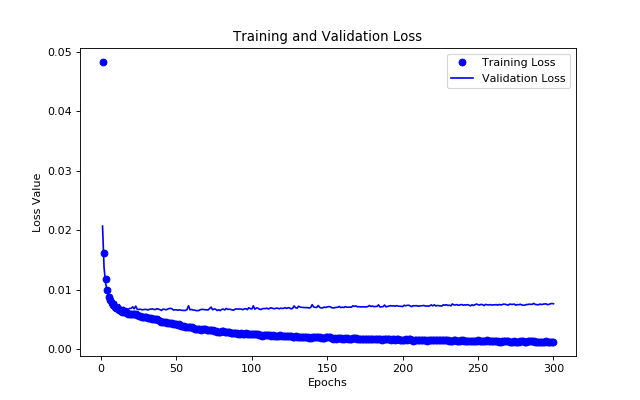}
\end{subfigure}
\begin{subfigure}{0.31\textwidth}
\centering
\footnotesize{Epochs -- 600 / Batch size - 4}\\
\includegraphics[width=0.99\linewidth]{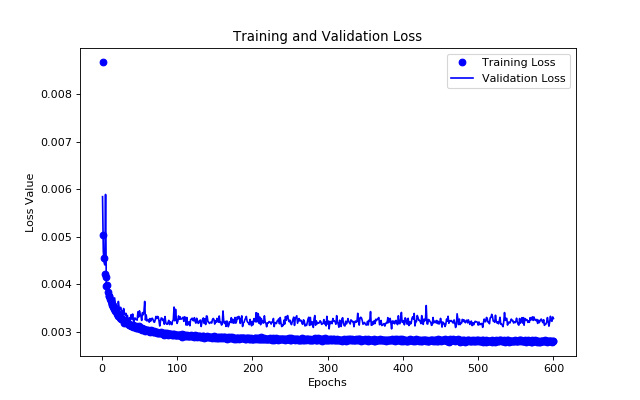}
\includegraphics[width=0.99\linewidth]{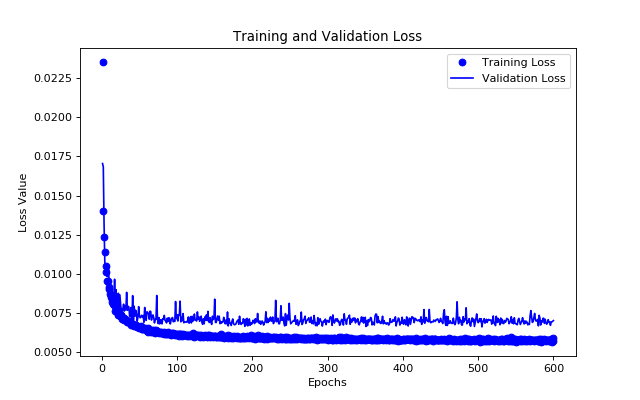}
\end{subfigure}
\begin{subfigure}{0.31\textwidth}
\centering
\footnotesize{Epochs -- 600 / Batch size - 16}\\
\includegraphics[width=0.99\linewidth]{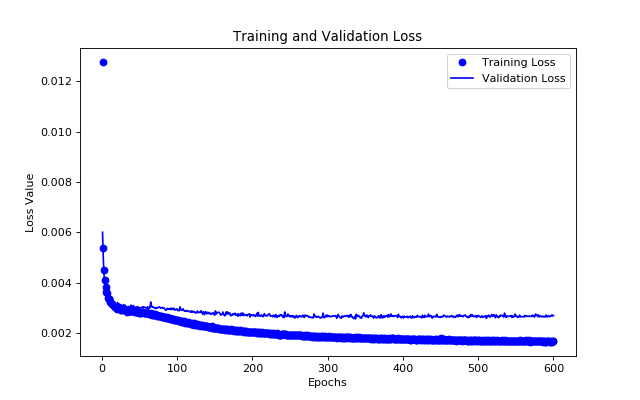}
\includegraphics[width=0.99\linewidth]{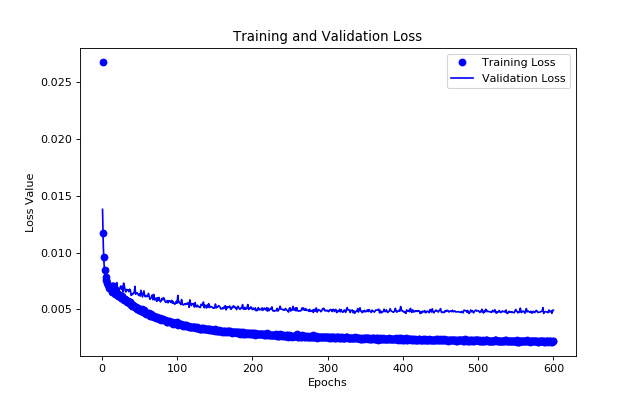}
\end{subfigure}
\begin{subfigure}{0.31\textwidth}
\centering
\footnotesize{Epochs -- 600 / Batch size - 64}\\
\includegraphics[width=0.99\linewidth]{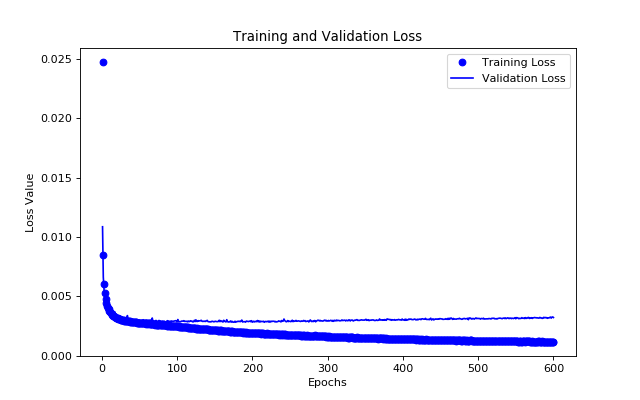}
\includegraphics[width=0.99\linewidth]{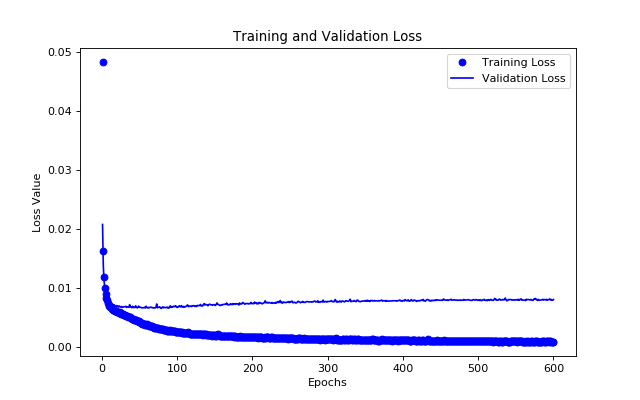}
\end{subfigure}

\caption{Steady-state case \eqref{eq:original00}.  Prediction's relative $l^2$ errors RMSE \eqref{rl2e} over 5000 training samples. First row: for effective permeability. Second row: for homogenized stiffness matrix.}
\label{PrErEp}
\end{figure}

Providing the steady-state form of the homogenized equation \eqref{rhomot}, 
we compute the macroscopic parameters (effective permeability tensor $\tilde{\kappa}$ \eqref{k2d} as well as coarse-scale homogenized stiffness matrix $\tilde{\bfa{A}}^n$ \eqref{comstiffh}) as in Subsection \ref{sec:method} and generate a data set of 6000 samples.  Then, 5000 training samples in such dataset are used to build deep neural networks (DNNs) for estimating those macroscopic parameters.
To show the convergence of loss function for the neural networks, Fig.~\ref{PrErEp} plots the relative $l^2$ error RMSE \eqref{rl2e} between the DNN-predictions and targets ($\tilde{\kappa}$ and $\tilde{\bfa{A}}^n$) in the training set (including the validation set)
corresponding to epoch number and batch size. 

Next, with the established deep neural networks, given 1000 testing pairs, we carry out the testing process for each macroscopic parameter. 
The minimum, maximum, mean, and one sample relative $l^2$ errors \eqref{Cer} between the DNN-predictions ($\tilde{\kappa}_{\tu{pred}}\,,$ $\tilde{\bfa{A}}^{n}_{\tu{pred}}$) and targets ($\tilde{\kappa}\,,$ $\tilde{\bfa{A}}^{n}$) are computed and recorded in Tables \ref{testres1-150}--\ref{testres1-600}.
As one can see from these tables, the predictions' accuracy improves as the number of epochs gets larger. We also take into account the increasing batch size, which leads to the growing learning rate as expected. 
With regard to accuracy, a very good neural network architecture is employed.


\begin{table}[H]
	\begin{center}
\begin{tabular}{cc}
		\begin{tabular}{|c|c|c|c|c|}
			\hline
			\parbox[c]{1.2cm}{Batch sizes} & min & max & mean & \parbox[c]{1.3cm}{one sample} \\
			\hline
			4 & 2.75 & 9.53 & 4.65 & 4.04\\
			8 & 2.96 & 8.47 & 4.63 & 4.19\\
			16 & 2.25 & 10.16 & 3.87 & 2.53\\
			32 & 2.03 & 6.96 & 3.14 & 2.21\\
			64 & 2.17 & 5.96 & 3.06 & 2.72\\
			\hline
		\end{tabular}
		&
		\begin{tabular}{|c|c|c|c|c|}
			\hline
			\parbox[c]{1.2cm}{Batch sizes} & min & max & mean & \parbox[c]{1.3cm}{one sample} \\
			\hline
			4 & 1.5 & 13.24 & 2.88 & 1.82 \\
			8 & 1.21 & 9.44 & 2.43 & 1.37 \\
			16 & 1.16 & 11.17 & 2.24 & 1.86 \\
			32 & 0.98 & 7.42 & 2.18 & 1.32 \\
			64 & 1.01 & 6.64 & 2.06 & 1.48 \\
			\hline
		\end{tabular}\\		
  (A) $e^{\tilde{\kappa}}_{l^2}$ (\%). & (B) $e^{\tilde{A}}_{l^2}$ (\%).
 	\end{tabular} 
	\end{center}
	\caption{Steady-state case \eqref{eq:original00}.  Prediction's relative $l^2$ errors \eqref{Cer} for 150 epochs with different batch sizes on 1000 testing samples: (A) for effective permeability (which is the same as time-dependent cases); (B) for homogenized stiffness matrix  (at the final Picard step).}
	\label{testres1-150}
\end{table}

\begin{table}[H]
	\begin{center}
\begin{tabular}{cc}
		\begin{tabular}{|c|c|c|c|c|}
			\hline
			\parbox[c]{1.2cm}{Batch sizes} & min & max & mean & \parbox[c]{1.3cm}{one sample} \\
			\hline
		4 & 2.95 & 10.08 & 4.63 & 3.91 \\
			8 & 3.03 & 9.86 & 4.98 & 4.47\\
			16 & 2.43 & 10.61 & 4.06 & 3.26 \\
			32 & 2.37 & 7.85 & 3.71 & 2.79\\
			64 & 1.89 & 7.93 & 3.09 & 2.21\\
			\hline
		\end{tabular}
		&
		\begin{tabular}{|c|c|c|c|c|}
			\hline
			\parbox[c]{1.2cm}{Batch sizes} & min & max & mean & \parbox[c]{1.3cm}{one sample} \\
			\hline
			4 & 1.45 & 12.62 & 2.53 & 2.12 \\
			8 & 1.24 & 7.63 & 2.31 & 1.51 \\
			16 & 1.26 & 11.61 & 2.35 & 1.26 \\
			32 & 1.15 & 6.05 & 2.06 & 1.27 \\
			64 & 0.85 & 8.66 & 1.96 & 1.06 \\
			\hline
		\end{tabular}\\		
  (A) $e^{\tilde{\kappa}}_{l^2}$ (\%). & (B) $e^{\tilde{A}}_{l^2}$ (\%).
 	\end{tabular} 
	\end{center}
	\caption{Steady-state case \eqref{eq:original00}.  Prediction's relative $l^2$ errors \eqref{Cer} for 300 epochs with different batch sizes on 1000 testing samples: (A) for effective permeability (which is the same as time-dependent cases); (B) for homogenized stiffness matrix  (at the final Picard step).}
	\label{testres1-300}
\end{table}

\begin{table}[H]
	\begin{center}
\begin{tabular}{cc}
		\begin{tabular}{|c|c|c|c|c|}
			\hline
			\parbox[c]{1.2cm}{Batch sizes} & min & max & mean & \parbox[c]{1.3cm}{one sample} \\
			\hline
			4 & 3.64 & 9.31 & 5.1 & 4.13\\
			8 & 3.03 & 9.06 & 4.5 & 3.34\\
			16 & 2.69 & 9.95 & 4.21 & 3.31\\
			32 & 2.11 & 8.39 & 3.52 & 3.03\\
			64 & 1.71 & 8.42 & 3.13 & 2.96 \\
			\hline
		\end{tabular}
		&
		\begin{tabular}{|c|c|c|c|c|}
			\hline
			\parbox[c]{1.2cm}{Batch sizes} & min & max & mean & \parbox[c]{1.3cm}{one sample} \\
			\hline
		4 & 0.88 & 8.66 & 1.96 & 2.44 \\
			8 & 1.25 & 11.95 & 2.19 & 1.56 \\
			16 & 1.06 & 8.95 & 2.22 & 1.58 \\
			32 & 1.03 & 9.23 & 2.08 & 1.39 \\
			64 & 0.78 & 9.71 & 1.95 & 1.46 \\
			\hline
		\end{tabular}\\		
  (A) $e^{\tilde{\kappa}}_{l^2}$ (\%). & (B) $e^{\tilde{A}}_{l^2}$ (\%).
 	\end{tabular} 
	\end{center}
	\caption{Steady-state case \eqref{eq:original00}.  Prediction's relative $l^2$ errors \eqref{Cer} for 600 epochs with different batch sizes on 1000 testing samples: (A) for effective permeability (which is the same as time-dependent cases); (B) for homogenized stiffness matrix  (at the final Picard step).}
	\label{testres1-600}
\end{table}

\newpage

\begin{figure}[H]
	\begin{center}
	  \includegraphics[width=0.85\linewidth]{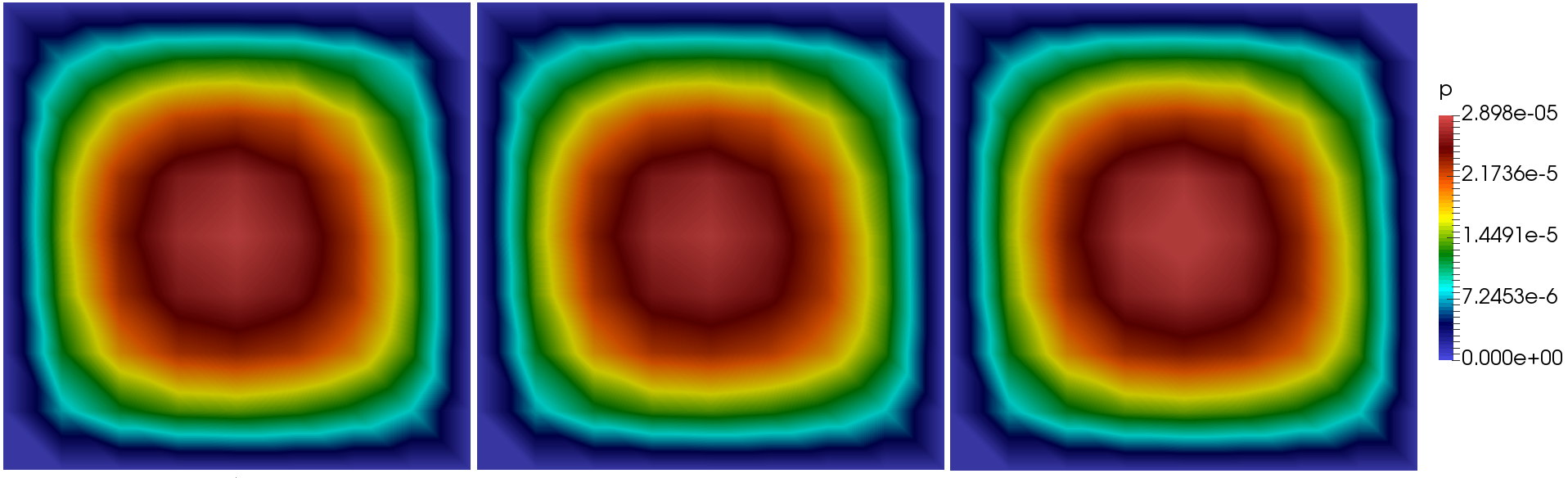}\\ (A) $p_c\,.$ \hspace{60pt} (B) $p^E_{c,\tu{pred}}\,.$ \hspace{50pt} (C) $p^A_{c,\tu{pred}}\,.$
	\end{center}
	\caption{Steady-state problem \eqref{eq:original00}.  Numerical results for solutions (using the first permeability field from Fig.~\ref{fig:kappa}): (A) homogenized solution;
	(B) solution obtained by prediction of effective permeability only; 
	(C) solution obtained by prediction of homogenized stiffness matrix alone.} 
	\label{fig:steady}
\end{figure}

Now, we present numerical results for solution to the homogenized equation of the steady-state problem \eqref{eq:original00}, using the first permeability field from Fig.~\ref{fig:kappa}. More specifically, Fig.~\ref{fig:steady} depicts the solution $p_c$ of \eqref{matrixf} obtained by numerical homogenization, solution $p^E_{c,\tu{pred}}$ to \eqref{apred} gained from DNN-prediction of effective permeability, and solution $p^A_{c,\tu{pred}}$ from \eqref{apred} attained by DNN-prediction of homogenized stiffness matrix.

Next, numerical experiments are completed on 1000 testing samples. That is, we compare the predicted solution with the homogenized solution $p_c$ \eqref{matrixf} in relative $L^2$ and $H^1$ errors \eqref{l2a}.  In particular, the comparisons between $p^E_{c,\tu{pred}}$ \eqref{apred} and $p_c$ \eqref{matrixf} are presented in Table \ref{tab2}, while the comparisons between $p^A_{c,\tu{pred}}$ \eqref{apred} and $p_c$ \eqref{matrixf} are shown in Table \ref{tab3}. In these tables, we depict mean, minimum, and maximum errors (over 1000 testing samples) as well as one sample errors (using the first permeability field in Fig.~\ref{fig:kappa}). It is observable from such results for solutions that DNN-prediction by effective permeability is a little bit more accurate than DNN-prediction by matrix. All predictions have good accuracy, and maximum errors are less than $10 \%$. The difference between mean errors and minimum errors are smaller than the difference between mean errors and maximum errors. This fact indicates that most of the numerical results (by accuracy) are closer to the minimum error. In summary, the DNN-prediction for steady-state Richards equation provides highly accurate solutions.

\begin{table}[H]
%
	\begin{center}
		\begin{tabular}{|c|c|c|}
			\hline
			 & $e^E_{L^2}$ (\%) & $e^E_{H^1}$ (\%)   \\
			\hline
			Mean error & 0.967 & 1.912 \\
			Minimum error & 0.305 & 0.843 \\
			Maximum error & 5.703 & 6.527 \\
			\hline
			One sample & 0.541 & 1.223 \\
			\hline
		\end{tabular}
		
	\end{center}
	\caption{Steady-state case \eqref{eq:original00}.  Numerical results for solutions (the case one sample uses the first permeability field from Fig.~\ref{fig:kappa}): errors \eqref{l2a} between homogenized solution and solution by predicted effective permeability.}
	\label{tab2}
\end{table}

\begin{table}[H]
%
	\begin{center}
		\begin{tabular}{|c|c|c|}
			\hline
			 & $e^A_{L^2}$ (\%) & $e^A_{H^1}$ (\%)   \\
			\hline
			Mean error & 1.238 & 2.431 \\
			Minimum error & 0.321 & 1.059 \\
			Maximum error & 5.676 & 8.034 \\
			\hline
			One sample & 0.676 & 1.778 \\
			\hline
		\end{tabular}
		
	\end{center}
	\caption{Steady-state case \eqref{eq:original00}.  Numerical results for solutions (the case one sample uses the first permeability field from Fig.~\ref{fig:kappa}): errors \eqref{l2a} between homogenized solution and solution by DNN-prediction of homogenized stiffness matrix.}
	\label{tab3}
\end{table}

\subsection{Experiment for time-dependent Richards equation}\label{tdc}

We investigate in this subsection the time-dependent Richards equation \eqref{eq:original00t} as a special case of \eqref{eq:original0}: find $p \in V$ such that 
\beq
\frac{\partial p(t,\bfa{x})}{\partial t} - \div\left ( \kappa(\bfa{x}) \, \frac{1}{1+|p|}\, \nabla p(t,\bfa{x}) \right )
= f(t,\bfa{x}) \  \textrm{in} \ (0,1] \times \Omega\,,
\eeq
where $f(t,\bfa{x}) = \cos(\pi x_1)\sin(\pi x_2)\,,$ with the initial condition $p(0,\bfa{x})= p_{0}=0$ over $\Omega\,,$ and in the company of the Dirichlet boundary condition $p(t,\bfa{x})=0$ on $(0,T] \times \partial \Omega\,.$  We set the terminal time $T=5\cdot 10^{-5} = S\tau$ and $S=20$ time steps, possessing the temporal step size $\tau=T/S = 25 \cdot 10^{-7}\,.$

Regarding the homogenized equation \eqref{rhomot} or \eqref{matrixf}, we calculate the macroscopic parameters (effective permeability tensor $\tilde{\kappa}$ \eqref{k2d}, homogenized stiffness matrix $\tilde{\bfa{A}}^n$ \eqref{comstiffh} (by data obtained from the $n$th Picard step), and right-hand side vector $\bfa{b}$ \eqref{comstiffh2} from \eqref{r1elhomoc}) as in Subsection \ref{sec:method} and make 6000 samples.  Using 5000 training samples, deep neural networks (DNNs) are created to approximate such macroscopic parameters.  

Then, with the constructed deep neural networks, their performance is checked on 1000 testing samples for each of the macroscopic parameters (effective permeability tensor, homogenized stiffness matrix, and right-hand side vector of \eqref{r1elhomoc}).  Note that the coarse-scale homogenized stiffness matrix $\tilde{\bfa{A}}^{n}_{\tu{pred}}$ \eqref{mpred} and right-hand side vector $\bfa{b}_{\tu{pred}}$ \eqref{bpred} of \eqref{r1elhomoc} are predicted for the last temporal step $(s+1)$th and at the last Picard iteration step $(n+1)$th, by data obtained from the $n$th Picard step.
Whereas, we predict the effective permeability field $\tilde{\kappa}_{\tu{pred}}$ \eqref{kpred} only once in Subsection \ref{ssc} within Tables \ref{testres1-150}--\ref{testres1-600}a, as there is not any time change. Tables \ref{testres3-150}--\ref{testres3-600} show the relative $l^2$ errors \eqref{Cer} between the targets ($\tilde{\bfa{A}}^{n}\,,$ $\bfa{b}$) and predictions ($\tilde{\bfa{A}}^{n}_{\tu{pred}}\,,$ $\bfa{b}_{\tu{pred}}$), 
at the last temporal step and final Picard step.  These Tables \ref{testres3-150}--\ref{testres3-600} demonstrate that the best results are obtained with 300 epochs and batch size 64, which will be used for prediction in the rest of this paper.  

\begin{figure}[H]
	\begin{center}
		\begin{minipage}[h]{0.45\linewidth}
			\center{\includegraphics[width=\linewidth]{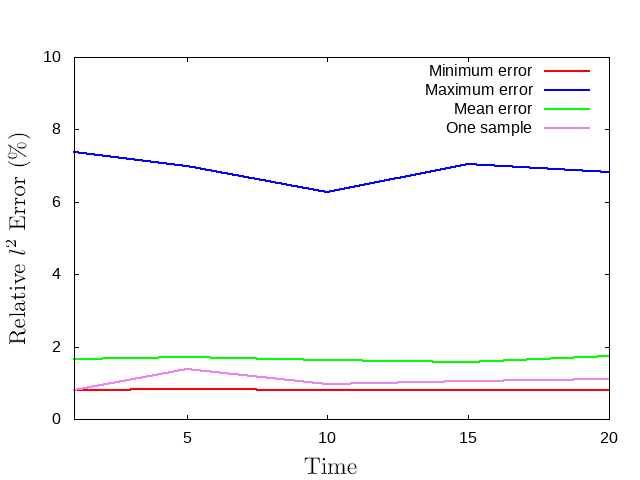}\\(A) $\tilde{\bfa{A}}^{n}_{\tu{pred}}\,.$}
		\end{minipage}
		\begin{minipage}[h]{0.45\linewidth}
			\center{\includegraphics[width=\linewidth]{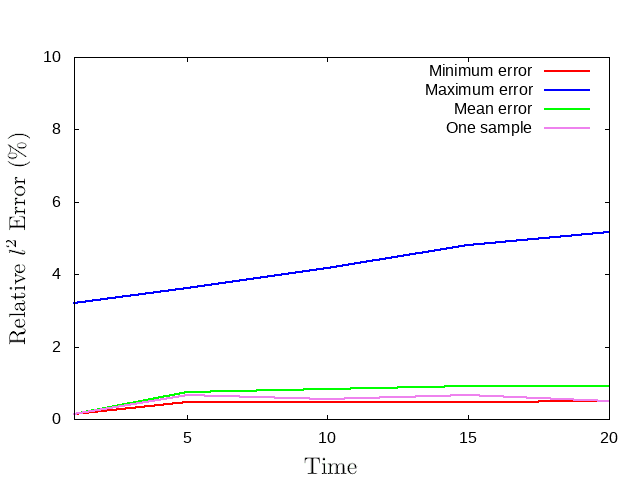}\\(B) $\bfa{b}_{\tu{pred}}\,.$}
		\end{minipage}
	\end{center}
	\caption{The time-dependent problem \eqref{eq:original00t}. Distribution of prediction's relative $l^2$ errors \eqref{Cer} at several time steps. Those are the min, max, and mean errors as well as one certain sample error: (A) for homogenized stiffness matrix; (B) for right-hand side vector of \eqref{r1elhomoc}.}
	\label{fig:preddm}
\end{figure}

For this time-dependent problem \eqref{eq:original00t}, the DNN-prediction by matrix and vector is also operated on 5 time steps (1st, 5th, 10th, 15th, 20th).  Over these temporal steps, Fig.~\ref{fig:preddm} presents the relative $l^2$ errors \eqref{Cer} between the targets ($\tilde{\bfa{A}}^{n}\,,$ $\bfa{b}$) and predictions ($\tilde{\bfa{A}}^{n}_{\tu{pred}}\,,$ $\bfa{b}_{\tu{pred}}$) over 1000 testing samples.  That errors are fairly small, and the largest errors are less than $8 \% \,.$

\begin{table}[H]
	\begin{center}
\begin{tabular}{cc}
		\begin{tabular}{|c|c|c|c|c|}
			\hline
			\parbox[c]{1.5cm}{Batch sizes} & min & max & mean & \parbox[c]{1.5cm}{one sample} \\
			\hline
			4 & 1.23 & 10.88 & 2.25 & 1.14 \\
			8 & 1.00 & 10.25 & 1.95 & 1.36 \\
			16 & 0.95 & 6.91 & 1.87 & 1.39 \\
			32 & 0.81 & 7.31 & 1.89 & 1.41 \\
			64 & 0.88 & 6.02 & 1.68 & 1.27 \\
			\hline
		\end{tabular}
  &
		\begin{tabular}{|c|c|c|c|c|}
			\hline
			\parbox[c]{1.5cm}{Batch sizes} & min & max & mean & \parbox[c]{1.5cm}{one sample} \\
			\hline
			4 & 0.65 & 5.19 & 1.05 & 1.19 \\
			8 & 0.61 & 5.95 & 1.05 & 0.77 \\
			16 & 0.59 & 4.36 & 1.15 & 1.11 \\
			32 & 0.56 & 4.78 & 0.99 & 0.74 \\
			64 & 0.77 & 6.7 & 1.27 & 1.16 \\
			\hline
		\end{tabular}\\
  (A) $e^{\tilde{A}}_{l^2}$ (\%). & (B) $e^{b}_{l^2}$ (\%).
 	\end{tabular} 
	\end{center}
\caption{The time-dependent problem \eqref{eq:original00t}.  Prediction's relative $l^2$ error \eqref{Cer} for 150 epochs with different batch sizes on 1000 testing samples: (A) for a homogenized stiffness matrix; (B) for a right-hand side vector of \eqref{r1elhomoc}.}
	\label{testres3-150}
\end{table}

\begin{table}[H]
	\begin{center}
\begin{tabular}{cc}
		\begin{tabular}{|c|c|c|c|c|}
			\hline
			\parbox[c]{1.5cm}{Batch sizes} & min & max & mean & \parbox[c]{1.5cm}{one sample} \\
			\hline
			4 & 1.21 & 10.93 & 2.20 & 2.00  \\
			8 & 1.15 & 9.54 & 2.24 & 1.92 \\
			16 & 0.87 & 9.11 & 1.68 & 1.14 \\
			32 & 0.89 & 6.84 & 1.74 & 1.29 \\
			64 & 0.81 & 7.16 & 1.63 & 1.03 \\
			\hline
		\end{tabular}
  &
		\begin{tabular}{|c|c|c|c|c|}
			\hline
			\parbox[c]{1.5cm}{Batch sizes} & min & max & mean & \parbox[c]{1.5cm}{one sample} \\
			\hline
			4 & 0.54 & 5.93 & 0.93 & 0.84 \\
			8 & 0.58 & 5.09 & 0.97 & 0.83 \\
			16 & 0.48 & 4.68 & 0.78 & 0.68 \\
			32 & 0.51 & 4.75 & 0.87 & 0.58\\
			64 & 0.47 & 5.16 & 0.91 & 0.47\\
			\hline
		\end{tabular}\\
  (A) $e^{\tilde{A}}_{l^2}$ (\%). & (B) $e^{b}_{l^2}$ (\%).
 	\end{tabular} 
	\end{center}
\caption{The time-dependent problem \eqref{eq:original00t}. Prediction's relative $l^2$ error \eqref{Cer} for 300 epochs with different batch sizes on 1000 testing samples: (A) for a homogenized stiffness matrix; (B) for a right-hand side vector of \eqref{r1elhomoc}.}
	\label{testres3-300}
\end{table}

\begin{table}[H]
	\begin{center}
\begin{tabular}{cc}
		\begin{tabular}{|c|c|c|c|c|}
			\hline
			\parbox[c]{1.5cm}{Batch sizes} & min & max & mean & \parbox[c]{1.5cm}{one sample} \\
			\hline
			4 & 1.61 & 10.9 & 2.46 & 2.14 \\
			8 & 0.94 & 13.87 & 1.76 & 1.27 \\
			16 & 1.01 & 7.08 & 1.95 & 1.29 \\
			32 & 0.81 & 6.81 & 1.65 & 0.88 \\
			64 & 0.78 & 8.55 & 1.69 & 0.98\\
			\hline
		\end{tabular}
  &
		\begin{tabular}{|c|c|c|c|c|}
			\hline
			\parbox[c]{1.5cm}{Batch sizes} & min & max & mean & \parbox[c]{1.5cm}{one sample} \\
			\hline
			4 & 0.61 & 3.22 & 1.01 & 0.93 \\
			8 & 0.55 & 2.95 & 0.93 & 0.66 \\
			16 & 0.48 & 4.26 & 0.74 & 0.56 \\
			32 & 0.51 & 4.01 & 0.84 & 0.55 \\
			64 & 0.51 & 4.41 & 0.82 & 0.56 \\
			\hline
		\end{tabular}\\
  (A) $e^{\tilde{A}}_{l^2}$ (\%). & (B) $e^{b}_{l^2}$ (\%).
 	\end{tabular} 
	\end{center}
\caption{The time-dependent problem \eqref{eq:original00t}.  Prediction's relative $l^2$ error \eqref{Cer} for 600 epochs with different batch sizes on 1000 testing samples: (A) for a homogenized stiffness matrix; (B) for a right-hand side vector of \eqref{r1elhomoc}.}
	\label{testres3-600}
\end{table}

\begin{figure}[H]
	\begin{center}
	  \includegraphics[width=\linewidth]{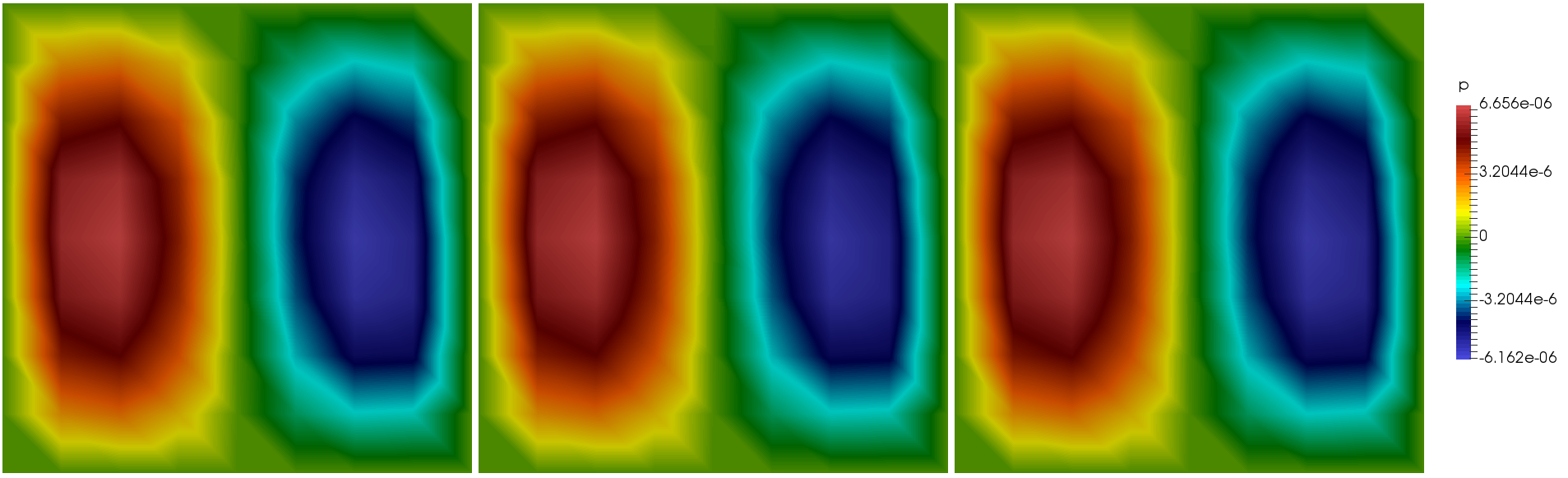}\\ (A) $p_c\,.$ \hspace{80pt} (B) $p^E_{c,\tu{pred}}\,.$ \hspace{80pt} (C) $p^A_{c,\tu{pred}}\,.$
	\end{center}
	\caption{The time-dependent problem \eqref{eq:original00t}.  Numerical results for solutions at the last temporal step (using the first permeability field from Fig.\ \ref{fig:kappa}): (A) homogenized solution; 
	(B) solution obtained by prediction of effective permeability;
	(C) solution obtained by prediction of the homogenized stiffness matrix
	and right-hand side vector
	of \eqref{r1elhomoc} at the last Picard step.}
	\label{fig:dynamic}
\end{figure}

Now, Fig.~\ref{fig:dynamic} shows numerical results for the solution to the homogenized equation of the time-dependent Richards equation \eqref{eq:original00t} at the last temporal step, using the first permeability field from Fig.~\ref{fig:kappa}: solution $p_c$ to \eqref{matrixf} is obtained by numerical homogenization,  solution $p^E_{c,\tu{pred}}$ of \eqref{apred} is achieved by DNN-prediction of effective permeability field, solution $p^A_{c,\tu{pred}}$ from \eqref{apred} is reached by DNN-prediction of the coarse-scale homogenized stiffness matrix together with right-hand side vector of \eqref{r1elhomoc} at the last Picard step.  

\begin{figure}[H]
	\begin{center}
		\begin{minipage}[h]{0.45\linewidth}
			\center{\includegraphics[width=0.92\linewidth]{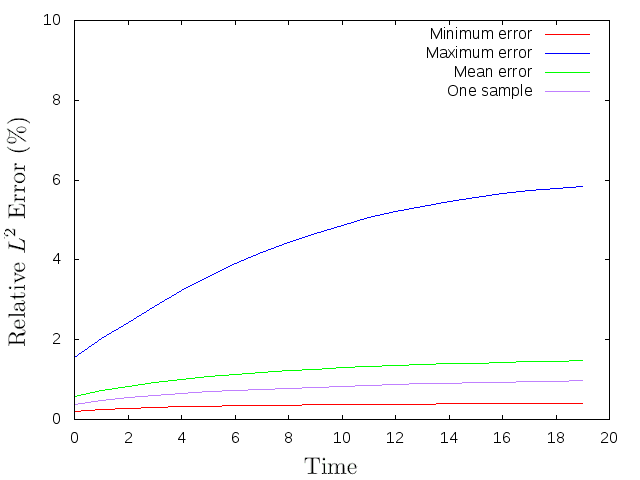}\\(A) Relative $L^2$ errors between $p_c$ and $p^E_{c,\tu{pred}}\,.$}
		\end{minipage}
		\begin{minipage}[h]{0.45\linewidth}
			\center{\includegraphics[width=0.92\linewidth]{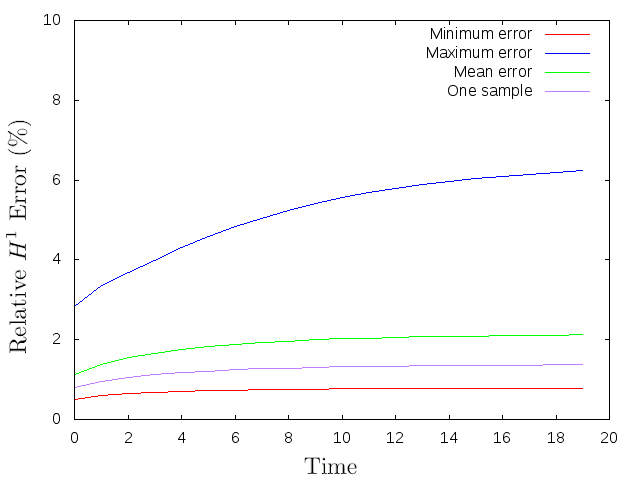}\\(B) Relative $H^1$ errors between $p_c$ and $p^E_{c,\tu{pred}}\,.$}
		\end{minipage}
	\end{center}
	\caption{The time-dependent problem \eqref{eq:original00t}.  Numerical results for solutions achieved by homogenization and by DNN-prediction of effective permeability (the case one sample uses the first permeability field from Fig.~\ref{fig:kappa}): (A) distribution of relative $L^2$ errors \eqref{l2a} by time; (B) distribution of relative $H^1$ errors \eqref{l2a} by time.}
	\label{graph1}
\end{figure}

\begin{figure}[H]
	\begin{center}
		\begin{minipage}[h]{0.45\linewidth}
			\center{\includegraphics[width=0.92\linewidth]{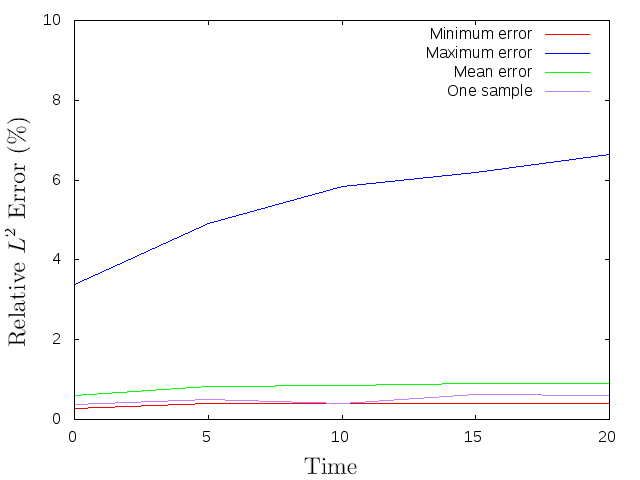}\\(A) Relative $L^2$ errors between $p_c$ and $p^A_{c,\tu{pred}}\,.$}
		\end{minipage}
		\begin{minipage}[h]{0.45\linewidth}
			\center{\includegraphics[width=0.92\linewidth]{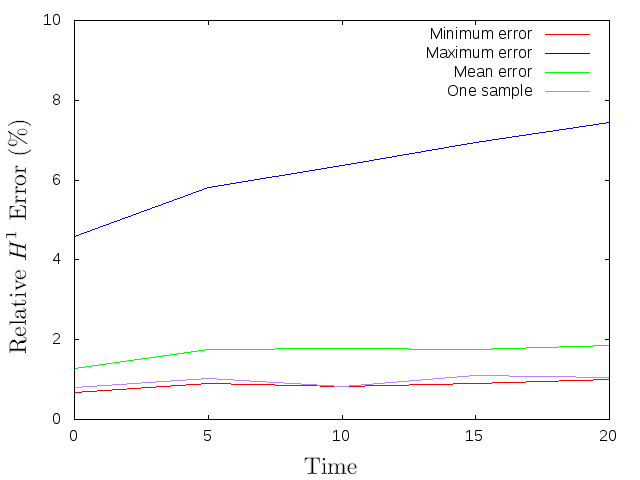}\\(B) Relative $H^1$ errors between $p_c$ and $p^A_{c,\tu{pred}}\,.$}
		\end{minipage}
	\end{center}
	\caption{The time-dependent problem \eqref{eq:original00t}.  Numerical results for solutions obtained by homogenization and by DNN-prediction of the coarse-scale homogenized stiffness matrix together with right-hand side vector of \eqref{r1elhomoc} (the case one sample uses the first permeability field from Fig.~\ref{fig:kappa}): (A) distribution of relative $L^2$ errors \eqref{l2a} by time; (B) distribution of relative $H^1$ errors \eqref{l2a} by time.}
	\label{graph4}
\end{figure}

Then, 1000 test samples are used in numerical experiments. In other words, we compare the homogenized solution $p_c$ \eqref{matrixf} with the predicted solution in terms of relative $L^2$ and $H^1$ errors \eqref{l2a}.  Tables \ref{tab5} and \ref{tab6} present the relative $L^2$ and $H^1$ errors between $p_c$ \eqref{matrixf} and $p^E_{c,\tu{pred}}$ \eqref{apred} as well as between $p_c$ \eqref{matrixf} and $p^A_{c,\tu{pred}}$ \eqref{apred}, respectively. We show mean, minimum, maximum, and one sample errors 
for the last Picard iteration at the last time step.  Figs.~\ref{graph1} and \ref{graph4} express the distribution of solutions' relative $L^2$ and $H^1$ errors over the chosen 5 time steps (1st, 5th, 10th, 15th, 20th).

It is clear from Tables \ref{tab5} and \ref{tab6} that the errors' behaviors are similar to the steady-state case \eqref{eq:original00}. The only difference is that for the solutions' mean errors, DNN-prediction by homogenized stiffness matrix and right-hand side vector of \eqref{r1elhomoc} becomes better than DNN-prediction by effective permeability, and this difference is also observable on Figs.~\ref{graph1} and \ref{graph4}.  Generally, however, both the approaches (prediction by effective permeability as well as by matrix and vector) provide solutions with small errors. For all cases, errors are smaller than $8 \%$. By accuracy, the mean error distribution is closer to the minimum error than to the maximum error. Note that for this time-dependent case, the approach with prediction of matrix and vector works better than the one with prediction of effective permeability because in applied problems, we need to know solution only at some important temporal points. It is more profitable to get only four or five crucial solutions than to make complex calculations over the entire temporal interval. In summary, the proposed approach has a good accuracy for the time-dependent Richards problem \eqref{eq:original00t}
. 

\begin{table}[H]
%
	\begin{center}
		\begin{tabular}{|c|c|c|}
			\hline
			 & $e^E_{L^2}$ (\%) & $e^E_{H^1}$ (\%)   \\
			\hline
			Mean error & 1.469 & 2.121 \\
			Minimum error & 0.403 & 0.781 \\
			Maximum error & 5.852 & 6.231 \\
			\hline
			One sample & 0.966 & 1.372 \\
			\hline
		\end{tabular}
	\end{center}
	\caption{The time-dependent problem \eqref{eq:original00t}.  Numerical results for solutions at the last time step (the case one sample uses the first permeability field from Fig.\ \ref{fig:kappa}): relative $L^2$ and $H^1$ errors \eqref{l2a} between homogenized solution $p_c$ and solution $p^E_{c,\tu{pred}}$ obtained by prediction of effective permeability at the last Picard step.}
	\label{tab5}
\end{table}

\begin{table}[H]
		
	\begin{center}
		\begin{tabular}{|c|c|c|}
			\hline
			 & $e^A_{L^2}$ (\%) & $e^A_{H^1}$ (\%)   \\
			\hline
			Mean error & 0.903 & 1.843 \\
			Minimum error & 0.399 & 1.011 \\
			Maximum error & 6.649 & 7.444 \\
			\hline
			One sample & 0.607 & 1.061 \\
			\hline
		\end{tabular}
		
	\end{center}
	\caption{The time-dependent problem \eqref{eq:original00t}.  Numerical results for solutions at the last time step (the case one sample uses the first permeability field from Fig.~\ref{fig:kappa}): relative $L^2$ and $H^1$ errors \eqref{l2a} between homogenized solution $p_c$ and solution $p^A_{c,\tu{pred}}$ obtained by prediction of the homogenized stiffness matrix together with right-hand side vector
	of \eqref{r1elhomoc} at the last Picard step.}
	\label{tab6}
\end{table}

\section{Conclusions}\label{sec:conclude}

In this paper, we establish a new coarse-grid approximation algorithm for the Richards equation as an unsaturated flow throughout heterogeneous media, utilizing numerical homogenization.
Our method takes advantage of deep neural networks (DNN) to swiftly and repeatedly compute macroscopic parameters. To be more precise, we train neural networks using a training set made up of stochastic permeability realizations and corresponding calculated macroscopic parameters (effective permeability, homogenized stiffness matrix, and right-hand side vector). Our suggested deep learning method creates maps between permeability fields and such macroscopic features, and the handling of Richards equation's nonlinearity is incorporated in the prediction of coarse-scale homogenized stiffness matrix, which is a novelty.  Finally, some good numerical results, especially for solutions, are presented to confirm the theory.

\bigskip

\noindent \textbf{Acknowledgements.} 
The reported study was funded by RFBR-VAST under the project number 21-51-54001 and supported by Duy Tan University under decision 5390/QD-DHDT.


\appendix

\section{Global convergence of Picard linearization procedure}\label{cp}
%

Consider Eq.~\eqref{eq:original0}, following \cite{rtt21, rpicardc,LPicard} (as well as thanks to J. Batista and A. Mazzucato, personal
communication, January 9, 2022), we will prove the global convergence of Picard linearization process (presented in Subsection \ref{pre}) for the system \eqref{r1elH} based on \eqref{r1ed}, with $i = 1,2$.  For simplicity, we remove from \eqref{r1elH} the subscripts $H$ and $(s+1)\,,$ where Picard iteration steps are applied.  Let $\| \cdot \|$ and $(\cdot , \cdot)$ respectively denote the $L^2(\Omega)$-norm and inner product in $L^2(\Omega)$.

Recall from Assumption \eqref{Coercivity} that the conductivity coefficient $\varkappa$ satisfies $ \underline{\varkappa} \leq \varkappa \leq \overline{\varkappa}\,,$ for some constants $\underline{\varkappa}, \overline{\varkappa} >0\,.$  Also, \eqref{ini} assumed that $p_0 = p(0,\bfa{x}) \in V = H_0^1(\Omega)\,.$  Without confusion, now $\varkappa(\bfa{x},p)$ is denoted by $\varkappa(p)\,.$  An additional assumption on $\varkappa(\bfa{x},\cdot)$ is the Lipschitz continuity as follows.
\begin{assumption}\label{ak}
There exists $L_{\varkappa} > 0$ (without any explicit dependence on $\bfa{x}$ and $t$) such that for all $\bfa{x} \in \Omega$ and $p_1,p_2 \in V\,,$
\begin{equation}\label{klip}
|\varkappa(p_1) - \varkappa(p_2)| \leq L_{\varkappa} |p_1 - p_2|\,.
\end{equation}
\end{assumption} 

By Poincar\'{e} inequality, there exists a constant $C_{\Omega}$ depending only on $\Omega$ such that 
\begin{equation}\label{poincare}
\|p^{n+1} - p_{s+1}\| \leq C_{\Omega}  \|\nabla (p^{n+1} - p_{s+1})\| \,,
\end{equation}
Also, assume that $p_{s+1}(t, \cdot) \in C_c^{\infty}(\Omega)$ so that $\|p_{s+1} \|_{\infty} \leq \hat{D} \|\nabla  p_{s+1} \|_{\infty}\,,$ 
where $\hat{D}$ is the distance between the two parallel hyperplanes bounding $\Omega\,.$  
Let $U$ be the exact solution of the equation at hand \eqref{eq:original0}.  Then, we obtain the following inequalities by applying the error estimate in \cite{GalerkinFEM} (Theorem 1.5):
\begin{align}\label{dsol}
\begin{split}
\frac{1}{\hat{D}} \|  p_{s+1} \|_{\infty} \leq\| \nabla p_{s+1} \|_{\infty}
& \leq \|\nabla(p_{s+1} - U(t_{s+1})) \|_{\infty} + \|\nabla U(t_{s+1})\|_{\infty}\\ 
& \leq \overline{C}(U)(H + \tau) + \|\nabla U(t_{s+1})\|_{\infty} = M\,,
\end{split}
\end{align}
for some constant $\overline{C}(U)$ depending on $U\,.$

Then in $V^H\,,$ we get the following equality by subtracting \eqref{r1ed} from \eqref{r1elH} and picking correct $v = p^{n+1} - p_{s+1}\,:$
\begin{align}\label{in1H}
\frac{1}{\tau} \|p^{n+1} - p_{s+1}\|^2 + a(p^{n+1},p^{n+1} - p_{s+1};p^n) = a(p_{s+1},p^{n+1} - p_{s+1};p_{s+1}) \,.
\end{align}
Equivalently,
\begin{align}\label{in1Hsub}
\begin{split}
\|p^{n+1} - p_{s+1}\|^2 + & \tau (\varkappa(p^n) \nabla(p^{n+1} - p_{s+1}), \nabla(p^{n+1} - p_{s+1})) \\
&=  -\tau ((\varkappa(p^n) - \varkappa(p_{s+1}))\nabla p_{s+1}, \nabla (p^{n+1} - p_{s+1})) \,.	
\end{split}
\end{align}	

Using Assumption \eqref{Coercivity} for the left-hand size of \eqref{in1Hsub}, then Cauchy-Schwarz inequality, Assumption \eqref{klip}, and the Young's inequality for the right-hand side of \eqref{in1Hsub}, we achieve
\begin{align}\label{in1alH}
\begin{split}
|| p^{n+1} - p_{s+1} ||^2 + \tau \underline{\varkappa} || \nabla (p^{n+1} - p_{s+1}) ||^2 \leq & \frac{\tau L^2_{\varkappa}}{2 \underline{\varkappa}} ||\nabla p_{s+1}||^2_{\infty} || p^n - p_{s+1}||^2 \\
& + \frac{\tau \underline{\varkappa}}{2} || \nabla(p^{n+1} - p_{s+1}) ||^2\,.
\end{split}
\end{align}
Employing \eqref{dsol} and \eqref{poincare}, we derive from \eqref{in1alH} that
\begin{align}
\left(1+ \frac{\tau \underline{\varkappa}}{2C^2_{\Omega}}\right) ||p^{n+1} - p_{s+1}||^2 \leq \frac{\tau L^2_{\varkappa} M^2}{2 \underline{\varkappa}}||p^n - p_{s+1}||^2\,,
\end{align}
and thus,
\begin{align}
||p^{n+1} - p_{s+1}||^2 \leq \frac{C L^2_{\varkappa}}{2\underline{\varkappa}} \, \frac{\tau}{1+ \dfrac{\tau \underline{\varkappa}}{2 C^2_{\Omega}}} (\tilde{C}(H+\tau) +1)^2 ||p^n - p_{s+1}||^2:= \lambda ||p^n - p_{s+1}||^2 \,,
\end{align}
for some positive constants $C, \tilde{C}\,.$
With sufficiently small $\tau$ and $H\,,$ the coefficient $\lambda$ will be less than $1\,,$ implying that the algorithm converges.  Specifically, $\lambda \to 0$ when $H \to 0$ and $\tau \to 0$ at the same time.

Similar proof holds for Eq.\ \eqref{eq:original00t} and its homogenized Eq.\ \eqref{rhomot} having discretization \eqref{r1elhomoc}.


\bibliographystyle{plain}
\bibliography{r1,r2}

\end{document}